\documentclass[12pt,a4paper,twoside]{article}
\usepackage[french] {babel}
\usepackage{amssymb,amsmath,amsfonts,graphics}
\usepackage{stmaryrd}
\usepackage[all]{xy}

\def\rondI{\newbox\boxx{\hbox{$1$}}\hskip-11pt\raisebox{0pt}{$\displaystyle{\textrm{\Large$\bigcirc$}}$}}
\def\rondII{\newbox\boxx{\hbox{$2$}}\hskip-11pt\raisebox{0pt}{$\displaystyle{\textrm{\Large$\bigcirc$}}$}}

\begin{document}

\newcommand{\hooklongrightarrow}{\lhook\joinrel\longrightarrow}
\begin{center}
\LARGE\textbf{Images directes III: F-isocristaux surconvergents }

\vskip20mm

Jean-Yves ETESSE
\footnote{(CNRS - Institut de Math\'ematique, Universit\'e de Rennes 1, Campus de Beaulieu - 35042 RENNES Cedex France)\\
E-mail : Jean-Yves.Etesse@univ-rennes1.fr}
\end{center}

 \vskip60mm

\noindent\textbf{Sommaire}\\

		\begin{enumerate}
		\item[] Introduction
		\item[1.] Frobenius
		\item[2.] Cas relevable
		\item[3.] Cas propre et lisse
		\item[4.] Cas fini \'etale
		\item[5.] Cas plongeable
		\end{enumerate}

\newpage

\noindent\textbf{R\'esum\'e}\\

Cet article est le troisi\`eme  d'une s\'erie de trois articles consacr\'es aux images directes d'isocristaux: ici nous consid\'erons des isocristaux surconvergents avec structure de Frobenius.\\

Pour un morphisme propre et lisse relevable nous \'etablissons la surconvergence des images directes, gr\^ace \`a des rel\`evements du Frobenius et au premier article. Ce r\'esultat r\'epond partiellement \`a une conjecture de Berthelot sur la surconvergence des images directes des $F$-isocristaux surconvergents par un morphisme propre et lisse.\\
  
\noindent\textbf{Abstract}\\

This article is the third one of a series of three articles devoted to direct images of isocrystals: here we consider overconvergent isocrystals with Frobenius structure.\\

For a liftable proper smooth morphism we establish the overconvergence of direct images, owing to the first article and the existence of lifts of Frobenius. This result partially answers a conjecture of Berthelot on the overconvergence of direct images of overconvergent $F$-isocrystals under a proper smooth morphism.\\


Mots cl\'es: alg\`ebres de Monsky-Washnitzer, sch\'emas formels, espaces rigides analytiques, $F$- isocristaux surconvergents, cohomologie rigide, images directes.\\


\newpage
\section*{ Introduction}

Cet article est le troisi\`eme d'une s\'erie de trois articles consacr\'es aux images directes d'isocristaux. Ici  on \'etudie les images directes des $F$-isocristaux dans le cadre \guillemotleft surconvergent\guillemotright. Par rapport \`a [Et 7] ([Et 5, chap II]) on ajoute donc une structure de Frobenius: tout le probl\`eme est d'obtenir de \guillemotleft bons\guillemotright \  rel\`evements du Frobenius.\\

Soient $\mathcal{V}$ un anneau de valuation discr\`ete complet, de corps r\'esiduel $k = \mathcal{V}/\mathfrak{m}$ de caract\'eristique $p > 0$, de corps des fractions $K$ de caract\'eristique $0$  et d'indice de ramification $e$, et $S$ un $k$-sch\'ema affine et lisse. 
Sur la base $S$, la question du rel\`evement du Frobenius est r\'esolue par le diagramme (1.2.4), qu'il s'agit ensuite de \guillemotleft tirer en haut\guillemotright \ par le morphisme propre et lisse $f:X\rightarrow S.$ \\

Dans le cas o\`u $f$ est relevable, on arrive \`a effectuer en parrall\`ele des rel\`evements du Frobenius et des bons choix de compactifications. D'o\`u le th\'eor\`eme de surconvergence des images directes dans le cas relevable [th\'eo (2.1)]: le point cl\'e est de montrer que le Frobenius, sur les images directes surconvergentes, est un isomorphisme; par le th\'eor\`eme de changement de base on est ramen\'e \`a la m\^eme propri\'et\'e dans le cas convergent vue dans [Et 8, (3.3.1)] ou [Et 5, chap III, (3.3.1)]. Dans le cas o\`u $f$ est projectif lisse relevable, ou $X$ intersection compl\`ete relative dans des espaces projectifs sur $S$, abord\'es au \S3, ou fini \'etale au \S4, on construit comme dans [Et 7, \S3] ou [Et 5, chap II, \S3] des foncteurs $R^{i}f_{rig\ast}$ sur la cat\'egorie des $F$-isocristaux surconvergents [th\'eor\`emes (3.1) et (4.1)].\\

Lorsque $f$ est propre et lisse et que de plus
\begin{enumerate}
\item[(i)] ou bien $f$ est g\'en\'eriquement projectif relevable, 
\item[(ii)] ou bien  $f$ est g\'en\'eriquement projectif et $X$ est intersection compl\`ete relative dans des espaces projectifs sur $S$,
\end{enumerate}
on utilise les r\'esultats de [Et 8]: on est ainsi amen\'e \`a supposer $k$ parfait, $e\leqslant p-1$ et \`a se restreindre \`a des $F$-isocristaux plats. Gr\^ace aux propri\'et\'es de descente \'etale des $F$-isocristaux surconvergents de [Et 4] et \`a la pleine fid\'elit\'e du foncteur d'oubli de la cat\'egorie surconvergente vers la cat\'egorie convergente, \'etablie par Kedlaya [Ked 2], on prouve encore la surconvergence des images directes [th\'eo (3.2)].\\

Dans le \S5, o\`u $f$ est suppos\'e seulement plongeable, on proc\`ede diff\'eremment: comme il n'y a plus de rel\`evements globaux du Frobenius comme pr\'ec\'edemment on utilise la fonctorialit\'e des images directes pour construire un morphisme de Frobenius; il reste alors \`a prouver que c'est un isomorphisme. Par un r\'esultat de Berthelot il suffit de voir que tel est le cas dans la cat\'egorie convergente: le th\'eor\`eme de changement de base et un r\'esultat de Bosch-G\"untzer-Remmert [B-G-R] nous ram\`ene \`a le v\'erifier aux points ferm\'es de $S$, pour lesquels l'assertion est connue [Et 8, (3.3.1.18)] ou [Et 5, chap III, (3.3.1.18)]. Pour $f$ plongeable on a ainsi des foncteurs images directes sur la cat\'egorie des $F$-isocristaux surconvergents [th\'eo (5.2)].

\section*{1. Frobenius}

\textbf{1.1.} On fixe dans ce paragraphe 1 un entier $a \in \mathbb{N}^{\ast}$ ; on pose $q = p^a$. Pour tout $k$-sch\'ema $S$ on notera $F_{S}$ le Frobenius de $S$ induit par la puissance $q$ sur le faisceau $\mathcal{O}_{S}$.\\

On fixe un rel\`evement $\sigma : \mathcal{V} \rightarrow \mathcal{V}$ de la puissance $q$ sur $k$ \`a la mani\`ere de [Et 4, I, 1.1].\\

Si $S$ est lisse sur $k$ et $e \leqslant p-1$, on notera  $F^{a}\mbox{-}\textrm{Isoc}^{\dag}(S/K)_{\textrm{plat}}$ la sous-cat\'egorie pleine de $F^{a}\mbox{-}\textrm{Isoc}^{\dag}(S/K)$ form\'e des objets dont l'image par le foncteur d'oubli

$$
F^{a}\mbox{-}\textrm{Isoc}^{\dag}(S/K) \longrightarrow F^{a}\mbox{-}\textrm{Isoc}(S/K) 
$$

\noindent est dans  $F^{a}\mbox{-}\textrm{Isoc}(S/K)_{\textrm{plat}}$ (cf [Et 8, 3.1] ou [Et 5, chap III, 3.1]).\\

\textbf{1.2.} Soit $S$ un $k$-sch\'ema affine et lisse. En utilisant les notations du [Et 6, th\'eo (3.1.3) (3)] ou [Et 5, chap I, th\'eo (3.4)(3)] il existe une $\mathcal{V}$-alg\`ebre lisse $A$ telle que $Spec\ A$ rel\`eve $S$ et un $\mathcal{V}$-morphisme fini $\psi$ relevant le Frobenius $F_{S}$, s'ins\'erant dans un diagramme commutatif \`a carr\'es cart\'esiens

$$
\begin{array}{c}
\xymatrix{
Spec\ B_{T} \ar[r] \ar[d]_{\psi_{T}}& Spec\ B\ \ar@{^{(}->}[r]^(.62){j_{Z'}} \ar[d]^{\psi} & Z' \ar[d]^{\overline{\psi}} \\
Spec\ A_{T} \ar[r]  & Spec\ A\ \ar@{^{(}->}[r]^(.62){j_{Z}}& Z
}
\end{array}
\leqno{(1.2.1)}
$$

\noindent o\`u les $j$ sont des immersions ouvertes, $\overline{\psi}$ est fini,$\psi_{T}$ est fini et plat, et $Z$ est un $\mathcal{V}$-sch\'ema propre, normal.\\

Soit $\hat{A}$ le s\'epar\'e compl\'et\'e $\mathfrak{m}$-adique de $A$ : c'est aussi le s\'epar\'e compl\'et\'e de $A_{T}$, et on a un isomorphisme [Et 6, th\'eo (3.1.3) (2) (i)] ou [Et 5, chap I, th\'eo (3.4)(2)(i)]

$$
B_{T} \otimes_{A_{T}} \hat{A} \simeq \hat{B} \simeq \hat{A} 
$$

\noindent tel que dans le diagramme commutatif \`a carr\'es cart\'esiens\\

$$
\begin{array}{c}
\xymatrix{
Spec\ \hat{A} \ar[r]^{\sim} \ar[rd]_{\varphi} \ar@/^2pc/[rr]^{\rho_{B}} & Spec\ \hat{B}_{T} \ar[r] \ar[d]^{\hat{\psi}_{T}} & Spec\ B\ \ar@{^{(}->}[r]^(.62){j_{Z'}} \ar[d]^{\psi} & Z' \ar[d]^{\overline{\psi}} \\
& Spec\ \hat{A}_{T} = Spec \hat{A} \ar[r]_(.62){\rho_{A}}  & Spec\ A\ \ar@{^{(}->}[r]_(.62){j_{Z}}& Z
}
\end{array}
\leqno{(1.2.2)}
$$

\noindent $\varphi$ est un rel\`evement de $F_{S}$.\\

Le morphisme diagonal $Spec\ \hat{A} \longrightarrow Spec\ \hat{A}\  \times_{\mathcal{V}}\  Spec\ \hat{A}$ est une immersion ferm\'ee, donc $Spec\ \hat{A}$ est isomorphe \`a son image sch\'ematique $\hat{\Delta}$ par ce morphisme. Consid\'erons l'image sch\'ematique de $\hat{\Delta}$ par le morphisme compos\'e

$$
Spec\ \hat{A}\ \times_{\mathcal{V}}\ Spec\ \hat{A} \displaystyle \mathop{\twoheadrightarrow}_{\rho=\rho_{_{A}}\times \rho_{B}}\ Spec\ A \times_{\mathcal{V}}\  Spec\ B \displaystyle \mathop{\hooklongrightarrow}_{j_{\mathcal{Z}} \times j_{\mathcal{Z}'} = j} \mathcal{Z} \times_{\mathcal{V}} \mathcal{Z}'\ ;
$$

\noindent notons $\Delta$ (resp. $\mathcal{Z}''$) l'image sch\'ematique de $\hat{\Delta}$  (resp. de $\Delta$) par $\rho_{A} \times \rho_{B}$ (resp. par $j_{Z} \times j_{Z'}$). L'immersion ouverte $j$ induit une immersion ouverte $j_{Z}'' : \Delta \hookrightarrow Z''$ [EGA I, (5.4.4)]. \\

Montrons que $\rho(\hat{\Delta}) = \Delta$. Quitte \`a d\'ecomposer la $\mathcal{V}$-alg\`ebre lisse (donc normale) $A$ en somme de ses composantes connexes, on peut supposer $A$ int\`egre, donc int\'egralement clos : ainsi $\hat{A}$ est int\'egralement clos [Et 6, prop (1.6) (4) (iv)] ou [Et 5, chap I, prop (1.6) (4) (iv)]. Soit $I$ l'id\'eal de $\hat{A} \otimes_{\mathcal{V}} \hat{A}$ d\'efinissant $\hat{\Delta} = Spec(\hat{A} \otimes_{\mathcal{V}} \hat{A}/I)$ : comme $\hat{A}$ est int\`egre, $I$ est un id\'eal premier. L'image de $A$ par $\rho$ est donc l'ensemble des id\'eaux premiers de $A \otimes_{\mathcal{V}} B$ contenant $\rho (I)$ : c'est donc $Spec(A \otimes_{\mathcal{V}} B/ \tilde{\rho}^{-1}(I))$ o\`u $\tilde{\rho} : A \otimes_{\mathcal{V}} B \rightarrow \hat{A} \otimes_{\mathcal{V}} \hat{A}$ induit $\rho$ ; comme c'est d\'ej\`a un sous-sch\'ema ferm\'e de $Spec\ A \times_{\mathcal{V}} Spec\ B$, il est \'egal \`a $\Delta$.\\

En remarquant que $\rho_{n} = \rho\ \textrm{mod}\  \mathfrak{m}^{n+1}$ est l'identit\'e, $\rho_{n}$ induit un isomorphisme

$$
(Spec\ \hat{A}\  \textrm{mod}\ \mathfrak{m}^{n+1}) \displaystyle \mathop{\longrightarrow}^{\sim} (\hat{\Delta} \ \textrm{mod}\ \mathfrak{m}^{n+1})  \displaystyle \mathop{\longrightarrow}^{\sim}_{\rho_{n}}
(\Delta\  \textrm{mod}\   \mathfrak{m}^{n+1}).
$$

Notons alors $\mathcal{S}, \mathcal{Z}, \mathcal{Z'}, \mathcal{Z''}$ les sch\'emas formels associ\'es respectivement \`a $Spec\ \hat{A}, \ Z, \ Z', \ Z''$. Ce qui pr\'ec\`ede fournit un diagramme commutatif \`a carr\'e cart\'esien\\

$$
\begin{array}{c}
\xymatrix{
& \mathcal{Z}\\
& \mathcal{Z}''  \ar@{^{(}->}[r] \ar[d] \ar[d]^{v_{\mathcal{Z}'}}  \ar[u]_{v_{\mathcal{Z}}} & \mathcal{Z} \times_{\mathcal{V}} \mathcal{Z}' \ar[dl]^{\textrm{proj}}  \ar[ul]_{\textrm{proj}}\\
\mathcal{S}   \ar@{^{(}->}[ur]_{j_{\mathcal{Z}''}}    \ar@{^{(}->}[uur]^{j_{\mathcal{Z}}}  \ar@{^{(}->}[r]_{j_{\mathcal{Z}'}}    \ar[d]_{F_{\mathcal{S}}:=\hat{\varphi}} &  \mathcal{Z}' \ar[d]^{\hat{\overline{\psi}}=:F_{\mathcal{Z}'}} &\\
\mathcal{S}  \ar@{^{(}->}[r]^{j_{\mathcal{Z}}}   & \mathcal{Z} & 
}
\end{array}
\leqno{(1.2.3)}  
$$

\noindent o\`u $\mathcal{Z}$ est propre sur $\mathcal{V}$ et normal et  $\hat{\overline{\psi}}$ est fini [Et 6, th\'eo (3.1.3)(3)] ou [Et 5, chap I, th\'eo (3.4) (3)]; $F_{\mathcal{S}}$ est un rel\`evement fini et plat du Frobenius $F_{S}$ ; $j_{\mathcal{Z}}, j_{\mathcal{Z'}}, j_{\mathcal{Z''}}$ sont des immersions ouvertes induites respectivement par $j_{Z}, j_{Z'}, j_{Z''}$ ; $i$ est l'immersion ferm\'ee induite par l'immersion ferm\'ee $\mathcal{Z}'' \hookrightarrow \mathcal{Z}\times_{\mathcal{V}}\mathcal{Z}'  ; \  v_{\mathcal{Z}}, v_{\mathcal{Z}'}$ sont des morphismes propres par composition de morphismes propres.\\

Avec les notations de [Et 7, (2.3.1)(2)] ou [Et 5, chap II, (2.3.1) (2)] soit $V_{\lambda} = Spm\ A_{\lambda}$: il existe $ \lambda_{0} > 1 $ tel que, pour $1< \lambda \leqslant \lambda_{0}$, $V_{\lambda}$ est lisse sur $K$ ; notons $W'_{\lambda} = F^{-1}_{\mathcal{Z}'_{K}}(V_{\lambda})$ et $W''_{\lambda} = v^{-1}_{\mathcal{Z}'_{K}}(W'_{\lambda})$. \\

\noindent Puisque $(V_{\lambda})_{\lambda}$ d\'ecrit un syst\`eme fondamental de voisinages stricts de $\mathcal{S}_{K}$ dans $\mathcal{Z}_{K}$, alors $(W'_{\lambda})_{\lambda}$ d\'ecrit un syst\`eme fondamental de voisinages stricts de $\mathcal{S}_{K}$ dans $\mathcal{Z}'_{K}$ [Et 7, prop (2.1.2)] ou [ Et 5, chap II, prop (2.1.2)]. Comme $v_{\mathcal{Z}'}$ est \'etale au voisinage de $\mathcal{S}$, il existe  $\lambda_{1} > 1$ tel que pour tout $\lambda$, $1 < \lambda \leqslant \lambda_{1} \leqslant\lambda_{0}$, on ait un isomorphisme $W''_{\lambda} \displaystyle \mathop{\rightarrow}^{\sim} W'_{\lambda}$ induit par $v_{\mathcal{Z}'}$ [B 3, (1.3.5)]. De m\^eme $v_{\mathcal{Z}}$ qui est \'etale au voisinage de $\mathcal{S}$, induit un isomorphisme entre un syst\`eme fondamental de voisinages stricts de $\mathcal{S}_{K}$ dans $\mathcal{Z}''_{K}$ et un syst\`eme fondamental de voisinages stricts de $\mathcal{S}_{K}$ dans $\mathcal{Z}_{K}$ : par composition il existe $\mu$, $1 < \mu \leqslant \lambda \leqslant \lambda_{1}$, et un morphisme fini $F_{\lambda \mu}$ rendant cart\'esien le carr\'e

$$
\begin{array}{c}
\xymatrix{
\mathcal{S}_{K} \ar@{^{(}->}[r] \ar[d]_{F_{\mathcal{S}_{K}}} & V_{\mu} \ar[d]^{F_{\lambda \mu}}\\
\mathcal{S}_{K} \ar@{^{(}->}[r] & V_{\lambda}
}
\end{array}
\leqno{(1.2.4)}
$$

\noindent o\`u les fl\`eches horizontales sont des immersions ouvertes et $V_{\lambda}$ est lisse sur $K$.

\vskip 6mm
\section*{2. Cas relevable}

\noindent \textbf{Th\'eor\`eme (2.1)}. \textit{Soient $S$ un $k$-sch\'ema lisse et s\'epar\'e et $f : X \rightarrow S$ un $k$-morphisme propre et lisse. On suppose qu'il existe un carr\'e cart\'esien de $\mathcal{V}$- sch\'emas formels}

$$
\begin{array}{c}
\xymatrix{
\mathcal{X} \ar@{^{(}->}[r]^{j_{\overline{\mathcal{X}}}} \ar[d]_{h} & \overline{\mathcal{X}} \ar[d]^{\overline{h}} &\\
\mathcal{S} \ar@{^{(}->}[r]_{j_{\overline{\mathcal{S}}}} & \overline{\mathcal{S}} & ,
}
\end{array}
\leqno{(2.1.1)}
$$

\noindent \textit{de r\'eduction mod $\mathfrak{m}$ \'egale \`a}

$$
\begin{array}{c}
\xymatrix{
X \ar@{^{(}->}[r]^{j_{\overline{X}}} \ar[d]_{f} & \overline{X} \ar[d]^{\overline{f}} &\\
S \ar@{^{(}->}[r]_{j_{\overline{S}}} & \overline{S} & ,
}
\end{array}
\leqno{(2.1.2)}
$$

\noindent \textit{o\`u $\overline{\mathcal{S}}$ est propre sur $\mathcal{V}$, $\overline{h}$ est propre, $h$ est propre et lisse et les $j$ sont des immersions ouvertes.}\\

\textit{Soit $E \in F^{a}\mbox{-}\textrm{Isoc}^{\dag}(X/K)$ et $\hat{E} = \mathcal{E}$ son image dans $F^{a}\mbox{-}\textrm{Isoc}(X/K)$. Alors, pour tout entier $i \geqslant 0$}

\begin{itemize}
\item[(1)] $E_{i} := R^{i} f_{\textrm{rig}^{\ast}}(X/ \overline{\mathcal{S}}, E) \in F^{a}\mbox{-}\textrm{Isoc}^{\dag}(S/K)$
\item[(2)] \textit{Soient $\hat{E}_{i} = j^{\ast}_{\overline{S}}(E_{i}), \mathcal{E}_{i} = R^{i} f_{\textrm{conv}^{\ast}} (X/ \mathcal{S}, \mathcal{E})$ et $\phi_{E_{i}} : F_{S}^{\ast}\ E_{i} \rightarrow E_{i}$, $\phi_{\mathcal{E}_{i}} : F_{S}^{\ast}\ \mathcal{E}_{i} \rightarrow \mathcal{E}_{i}$ les isomorphismes de Frobenius. Le diagramme commutatif d'isomorphismes ci-dessous d\'efinit $\phi_{\hat{E}_{i}}$ et permet les identifications canoniques} 
\end{itemize}

$$
j^{\ast}_{\overline{S}} (\phi_{E_{i}}) = \phi_{\hat{E}_{i}} = \phi_{\mathcal{E}_{i}}
$$

$$ 
\xymatrix{
j^{\ast}_{\overline{S}}\ F^{\ast}_{S}\ R^{i} f_{\textrm{rig} \ast}(X/ \overline{\mathcal{S}}, E) \ar[r]^{\sim}_(.53){j^{\ast}_{\overline{S}} (\phi_{E_{i}})} \ar[d]^{\simeq} & j^{\ast}_{\overline{S}}\ R^{i} f_{\textrm{rig}\ast}(X/ \overline{\mathcal{S}}, E) \ar@{=}[d]\\
F^{\ast}_{S}\ j^{\ast}_{\overline{S}}\  R^{i} f_{\textrm{rig} \ast}(X/ \overline{\mathcal{S}}, E) \ar[r]^{\sim}_{\phi_{\hat{E}_{i}}}   \ar[d]^{\simeq} & j^{\ast}_{\overline{S}}\  R^{i} f_{\textrm{rig}\ast}(X/ \overline{\mathcal{S}}, E) \ar[d]^{\simeq}\\
F^{\ast}_{S}\  R^{i} f_{\textrm{conv}^{\ast}}(X/ \mathcal{S}, \mathcal{E}) \ar[r]^{\sim}_{\phi_{\mathcal{E}_{i}}}  & R^{i}  f_{\textrm{conv}^{\ast}}(X/ \mathcal{S}, \mathcal{E}) .
}
$$

\vskip 3mm
\noindent \textit{D\'emonstration}. L'image inverse $F^{\ast}_{\sigma}\ E_{i}$ par Frobenius s'obtient [B 3, (2.3.7)] en appliquant le foncteur de changement de base
$$
\sigma^{\ast} : \textrm{Isoc}^{\dag}(S/K) \longrightarrow \textrm{Isoc}^{\dag}(S^{(q)}/K)
$$
\noindent puis le foncteur image inverse par le Frobenius $F_{S/k} : S \rightarrow S^{(q)}$.\\

\noindent Comme $\sigma$ est fix\'e on notera $F^{\ast}_{\sigma}\ E = F^{\ast}_{S} \ E$. Il nous reste donc \`a d\'efinir l'isomorphisme  de Frobenius $\phi_{E_{i}}$ de $E_{i}$.\\

Quitte \`a d\'ecomposer $S$ en somme de ses composantes connexes il suffit de d\'efinir $\phi_{E_{i}}$ sur chacune de ces composantes connexes. Soit $S_{\alpha}$ un ouvert affine d'une composante connexe $S_{0}$ de $S$ : comme le foncteur
$$
F^{a}\mbox{-} \textrm{Isoc}^{\dag}(S_{0}/K) \longrightarrow F^{a}\mbox{-} \textrm{Isoc}^{\dag}(S_{\alpha}/K)
$$
\noindent est pleinement fid\`ele [Et 4, th\'eo 4], il suffit de d\'efinir $\phi_{E_{i}}$ sur $S_{\alpha}$.\\

Soit $j_{s_{\alpha}} : S_{\alpha} = Spec\ A_{0} \hookrightarrow S$ l'immersion ouverte et $A$ une $\mathcal{V}$-alg\`ebre lisse relevant $A_{0}$. D'apr\`es (1.2.3) on a un diagramme commutatif de $\mathcal{V}$-sch\'emas formels de type fini, \`a carr\'e cart\'esien\\

$$
\begin{array}{c}
\xymatrix{
& & & &  \overline{\mathcal{S}}_{\alpha}\\
& & & \overline{\mathcal{S}}''_{\alpha}   \ar[ur]_{v_{\overline{\mathcal{S}}_{\alpha}}}  \ar[d]^{v_{\overline{\mathcal{S}}'_{\alpha}}} & &\\
S_{\alpha}  \ar@{^{(}->}[rrr]_{j_{\overline{\mathcal{S}}'_{\alpha}}}    \ar@{^{(}->}[urrr]_{j_{\overline{\mathcal{S}}''_{\alpha}}}      \ar@{^{(}->}[uurrrr]^{j_{\overline{\mathcal{S}}_{\alpha}}} \ar[d]_{F_{\alpha}} &  & &  \overline{\mathcal{S}}'_{\alpha} \ar[d]^{\overline{F}_{\alpha}}    &  \\
\mathcal{S}_{\alpha}  \ar@{^{(}->}[rrr]^{j_{\overline{\mathcal{S}}_{\alpha}}}  & & & \overline{\mathcal{S}}_{\alpha} & 
}
\end{array}
\leqno{(2.1.3)}  
$$

\noindent o\`u $\mathcal{S}_{\alpha} = Spf \hat{A}$, $\overline{\mathcal{S}}_{\alpha}$ est propre sur $\mathcal{V}$, $v_{\overline{\mathcal{S}}'_{\alpha}}$ et $v_{\overline{\mathcal{S}}_{\alpha}}$ sont propres, $F_{\alpha}$ est un rel\`evement fini et plat du Frobenius $F_{S_{\alpha}}$ de $S_{\alpha}$, $\overline{F}_{\alpha}$ est fini et les $j$ sont des immersions ouvertes. Notons $j_{\overline{\mathcal{T}}_{\alpha}} : \mathcal{T}_{\alpha} := \mathcal{S}_{\alpha} \times_{\mathcal{V}} \mathcal{S} \longrightarrow \overline{\mathcal{T}}_{\alpha} : = \overline{\mathcal{S}}_{\alpha} \times_{\mathcal{V}} \overline{\mathcal{S}}$ l'immersion ouverte et $\overline{u}_{\alpha} : \overline{\mathcal{T}}_{\alpha} \longrightarrow \overline{\mathcal{S}}_{\alpha}, \overline{v}_{\alpha} : \overline{\mathcal{T}}_{\alpha} \longrightarrow \overline{\mathcal{S}}$ les projections ; 
soient $\overline{T}_{\alpha}$ (resp. $T_{\alpha}$) la r\'eduction de $\overline{\mathcal{T}_{\alpha}}$ (resp. $\mathcal{T}_{\alpha}$) mod $\mathfrak{m}$ et $\tilde{S}_{\alpha}$ l'image sch\'ematique de $S_{\alpha}$ plong\'e diagonalement dans $\overline{T}_{\alpha}$ :

$$
\xymatrix{
S_{\alpha} \ar@{^{(}->}[r]_{j_{\tilde{S}_{\alpha}}} & \tilde{S}_{\alpha}  \ar@{^{(}->}[r]^{i_{\overline{T}_\alpha}}  \ar@/_1pc/[rr]_{i_{\tilde{S}_{\alpha}}} & \overline{T}_{\alpha} \ar@{^{(}->}[r]^{i_{\overline{\mathcal{T}}_\alpha}} & \overline{\mathcal{T}}_{\alpha}
}
$$

\noindent $j_{\tilde{S}_{\alpha}}$ est une immersion ouverte et les $i$ des immersions ferm\'ees. On a alors un diagramme commutatif \`a carr\'es verticaux cart\'esiens\\

$$
\begin{array}{c}
 \shorthandoff{;:!?}
 \xymatrix@!0 @R=1cm @C=2cm{
&&X\ar@{.>}[dd]^{f}\  \ar @{^{(}->}[rr]^{j_{\overline{X}}}&&\overline{X}  \ar@{.>}[dd]^{\overline{f}}\  \ar@{^{(}->}[rr]^{i_{\overline{X}}}&&\overline{\mathcal{X} }\ar[dd]^{\overline{h}}\\
&&&&&&\\
&&S\  \ar@{^{(}.>}[rr]^{j_{\overline{S}}}  && \overline{S}\  \ar@{^{(}.>}[rr]^{i_{\overline{S}}}&& \overline{\mathcal{S}}\\
&X_{\alpha}\ar@{.>}[dd]^{f_{\alpha}} \ar@{^{(}->}[uuur]^{j_{X_{\alpha}}} \ar@{^{(}->}[rr]  &&\overline{X}\ar@{.>}[dd]_{\overline{f}} \  \ar@{^{(}->}[rr]^(.7){i_{\overline{X}}} \ar@{=}[uuur]&& \overline{\mathcal{X}} \ar[dd]^{\overline{h}} \ar@{=}[uuur]\\
&&&&&&\\
& S_{\alpha}\  \ar@{^{(}.>}[rr]^{\overline{j}_{\alpha}} \ar@{.>}[uuur]_{j_{S_{\alpha}}} &&\overline{S} \  \ar@{^{(}.>}[rr] ^{i_{\overline{S}}} \ar@{.>}[uuur]_{id}&& \overline{\mathcal{S}} \ar@{=}[uuur]_{id}&\\
X_{\alpha} \ar@{^{(}->}[r] \ar[dd]_{f_{\alpha}}\   \ar@{=}[uuur]^{id}& \tilde{X}_{\alpha}\ar@{^{(}->}[r] \ar[dd]&\overline{Y}_{\alpha}\ar[dd]\  \ar@{^{(}->}[rr] \ar[uuur]&& \overline{\mathcal{Y}}_{\alpha}\ar[dd]^{\overline{h}_{\alpha}} \ar[uuur]&&\\
&&&&&&\\
S_{\alpha} \  \ar@{^{(}->}[r]_{j_{\tilde{S}_{\alpha}}}  \ar@{.>}[uuur]_{id}& \tilde{S}_{\alpha}\ar@{^{(}->}[r]&\overline{T}_{\alpha} \ar@{^{(}->}[rr] \  \ar@{.>}[uuur]&&\ \overline{\mathcal{T}}_{\alpha}\ . \ar[uuur]_{\overline{v}_{\alpha}}&&
}
\end{array}
\leqno{(2.1.4)}
$$

Ainsi on a une suite d'isomorphismes

$$
j^{\ast}_{S_{\alpha}}\ R^i f_{\textrm{rig}\ast}(X/\overline{\mathcal{S}}, E) \tilde{\longrightarrow} R^i f_{\alpha \textrm{rig}\ast} (X_{\alpha}/ \overline{\mathcal{S}}, j^{\ast}_{X_{\alpha}} E)\ [\textrm{[Et 7, th\'eo (3.4.4)] ou [Et 5, chap II, th\'eo(3.4.4)]}]
$$

\noindent (2.1.5) $\qquad \simeq \overline{v}^{\ast}_{\alpha_{K}} R^i f_{\alpha \textrm{rig}\ast}(X_{\alpha}/ \overline{\mathcal{S}}, j^{\ast}_{X_{\alpha}} E)\ [\textrm{B 3, (2.3.6), (2.3.2) (iv)}]$\\

$\qquad \qquad \tilde{\rightarrow} R^i f_{\alpha \textrm{rig}\ast}(X_{\alpha}/ \overline{\mathcal{T}_{\alpha}}, j^{\ast}_{X_{\alpha}} E)\ [\textrm{[Et 7, th\'eo (3.4.4)]}].$\\

\noindent On est donc ramen\'e \`a construire un isomorphisme de Frobenius sur

$$
R^i f_{\alpha \textrm{rig}\ast}(X_{\alpha}/ \overline{\mathcal{T}_{\alpha}}, j^{\ast}_{X_{\alpha}} E) \in \textrm{Isoc}^{\dag}(S_{\alpha}/K).
$$

\noindent A partir du carr\'e commutatif

$$
\xymatrix{
\mathcal{T}_{\alpha} \ar@{^{(}->}[rr]^(.4){j_{\overline{\mathcal{T}}'_{\alpha}}} \ar[d]_{F_{\mathcal{T}_{\alpha}} := F_{\alpha} \times 1} && \overline{\mathcal{T}}'_{\alpha} : = \overline{\mathcal{S}}'_{\alpha} \times_{\mathcal{V}} \overline{\mathcal{S}} \ar[d]^{\overline{F}_{\alpha} \times 1 =: F_{\overline{\mathcal{T}}'_{\alpha}}} &\\
\mathcal{T}_{\alpha} \ar@{^{(}->}[rr]_(.4){j_{\overline{\mathcal{T}}_{\alpha}}}  &&   \overline{\mathcal{T}}_{\alpha} = \overline{\mathcal{S}}_{\alpha} \times_{\mathcal{V}} \overline{\mathcal{S}} & ,
}
$$

\noindent d\'eduit de (2.1.3), et du carr\'e cart\'esien

$$
\xymatrix{
\mathcal{Y}_{\alpha} \ar@{^{(}->}[r] \ar[d]_{h_{\alpha}} & \overline{\mathcal{Y}_{\alpha}} \ar[d]^{\overline{h}_{\alpha}} &\\
\mathcal{T}_{\alpha} \ar@{^{(}->}[r]_{j_{\overline{\mathcal{T}}_{\alpha}}} & \overline{\mathcal{T}}_{\alpha} & ,
}
$$

\noindent o\`u les fl\`eches $j$ sont des immersions ouvertes, on forme les diagrammes commutatifs \`a carr\'es cart\'esiens\\

$$
\begin{array}{c}
\xymatrix{
& X_{\alpha} \ar@{.>}[dd]^(.3){f_{\alpha}} |\hole \ar[rr] & & Y_{\alpha} \ar@{.>}[dd] |\hole \ar@{^{(}->}[rr] & & \mathcal{Y}_{\alpha} \ar@{.>}[dd]^(.3){h_{\alpha}}  \ar@{^{(}->}[rr] && \overline{\mathcal{Y}}_{\alpha}\ar[dd]^{\overline{h}_{\alpha}} \\
X^{(q/S_{\alpha})}_{\alpha} \ar[rr] \ar[dd]_{f^{(q)}_{\alpha}} \ar[ur] & & Y'_{\alpha}\ar@{^{(}->}[rr] \ar[dd] \ar[ur]  & & \mathcal{Y}'_{\alpha} \ar[dd]  \ar[ur]  \ar@{^{(}->}[rr]  && \overline{\mathcal{Y}}'_{\alpha} \ar[ur] \ar[dd] \\
& S_{\alpha} \ar@{.>}[rr]^(.7){i_{\alpha}} |\hole & & T_{\alpha}\ar@{^{(}.>}[rr]^(.4){i_{\mathcal{T}_{\alpha}}}  |\hole & &  \mathcal{T}_{\alpha}  \ar@{^{(}.>}[rr]^(.4){j_{\overline{\mathcal{T}}_{\alpha}}}  |\hole && \overline{\mathcal{T}_{\alpha}}\\
S_{\alpha} \ar[rr]_{i_{\alpha}} \ar@{.>}[ur]^{F_{S_{\alpha}}} & & T_{\alpha} \ar@{^{(}->}[rr]_{i_{\mathcal{T}_{\alpha}}} \ar@{.>}[ur]_{F_{S_{\alpha}} \times 1} & & \mathcal{T}_{\alpha} \ar@{^{(}->}[rr]_{j_{\overline{\mathcal{T}}'_{\alpha}}} \ar@{.>}[ur]_{F_{\mathcal{T}_{\alpha}}} && \overline{\mathcal{T}}'_{\alpha} \ar[ur]_{F_{\overline{\mathcal{T}}'_{\alpha}}}
}
\end{array}
\leqno{(2.1.6)}
$$

\noindent (o\`u $i_{\alpha}$ et $i_{\mathcal{T}_{\alpha}}$ sont des immersions ferm\'ees), et \\

$$
\begin{array}{c}
\xymatrix{
& X_{\alpha} \ar@{.>}[dd]^(.3){f_{\alpha}} |\hole \ar@{^{(}->}[rr]^{j_{\tilde{X}_{\alpha}}} & & \tilde{X}_{\alpha} \ar@{.>}[dd]^(.3){\tilde{f}_{\alpha}} |\hole \ar[rr]& & \overline{\mathcal{Y}}_{\alpha} \ar[dd]^{\overline{h}_{\alpha}} \\
X^{(q/S_{\alpha})}_{\alpha} \ar@{^{(}->}[rr] \ar[dd]_{f^{(q)}_{\alpha}} \ar[ur]^{\pi_{X_{\alpha}}} & & \tilde{X}'_{\alpha} \ar[rr] \ar[dd]_(.3){\tilde{f}'_{\alpha}} \ar[ur]  & & \overline{\mathcal{Y}}'_{\alpha} \ar[dd]^(.3){\overline{h}'_{\alpha}} \ar[ur] \\
& S_{\alpha} \ar@{^{(}.>}[rr]^(.7){j_{\tilde{S}_{\alpha}}} |\hole & &\tilde{S}'_{\alpha}\ar@{^{(}.>}[rr]^(.4){i_{\tilde{S}_{\alpha}}}  |\hole & &  \overline{\mathcal{T}}'_{\alpha}  \\
S_{\alpha} \ar@{^{(}->}[rr]_{j_{\tilde{S}'_{\alpha}}} \ar@{.>}[ur]^{F_{S_{\alpha}}} & & \tilde{S}'_{\alpha} \ar@{^{(}->}[rr]_{i_{\tilde{S}'_{\alpha}}} \ar@{.>}[ur]_{\tilde{F}_{\alpha}} & & \  \overline{\mathcal{T}}'_{\alpha}\ . \ar[ur]_{F_{\overline{\mathcal{T}}'_{\alpha}}} 
}
\end{array}
\leqno{(2.1.7)}
$$

D'apr\`es [Et 7, th\'eo (3.4.4)] ou [Et 5, chap II, th\'eo (3.4.4)] on a un isomorphisme\\

\noindent (2.1.8)  $\qquad  F^{\ast}_{S_{\alpha}} R^i f_{\alpha\  \textrm{rig}^{\ast}} (X_{\alpha}/\overline{\mathcal{T}_{\alpha}}, j^{\ast}_{X_{\alpha}} E) \tilde{\rightarrow} R^i f^{(q)}_{\alpha\ \textrm{rig}^{\ast}} (X^{(q/S_{\alpha})}_{\alpha} / \overline{\mathcal{T}}'_{\alpha}, \pi^{\ast}_{X_{\alpha}}\ j^{\ast}_{X_{\alpha}} E).$ \\

Notons $v_{\overline{\mathcal{T}}_{\alpha}} := v_{\overline{\mathcal{S}}_{\alpha}} \times 1_{\overline{\mathcal{S}}} : \overline{\mathcal{T}}''_{\alpha} = \overline{\mathcal{S}}''_{\alpha} \times_{\mathcal{V}}    \overline{\mathcal{S}} \rightarrow  \overline{\mathcal{T}}_{\alpha} = \overline{\mathcal{S}}_{\alpha} \times_{\mathcal{V}} \overline{\mathcal{S}}$ et $v_{\overline{\mathcal{T}}'_{\alpha}} := v_{\overline{\mathcal{S}}'_{\alpha}} \times 1_{\overline{\mathcal{S}}} : \overline{\mathcal{T}}''_{\alpha} = \overline{\mathcal{S}}''_{\alpha} \times_{\mathcal{V}}  \overline{\mathcal{S}} \rightarrow \overline{\mathcal{T}}'_{\alpha} = \overline{\mathcal{S}}'_{\alpha} \times_{\mathcal{V}} \overline{\mathcal{S}}$ ; de (2.1.3) on d\'eduit un diagramme commutatif\\

$$
\xymatrix{
& & & & \overline{\mathcal{T}}_{\alpha}\\
S_{\alpha} \ar[r] &  \mathcal{T}_{\alpha} \ar@{^{(}->}[rr]^(.7){j_{\overline{\mathcal{T}}^{''}_{\alpha}}} \ar@{^{(}->}[rrrd]_{j_{\overline{\mathcal{T}}^{'}_{\alpha}}}   \ar@{^{(}->}[urrr]^{j_{\overline{\mathcal{T}}_{\alpha}}}  && \overline{\mathcal{T}}^{''}_{\alpha} \ar[rd]^{v_{\overline{\mathcal{T}}'_{\alpha}}}  \ar[ur]_{v_{\overline{\mathcal{T}}_{\alpha}}} \\
&&& &  \overline{\mathcal{T}}'_{\alpha}
}
$$

\noindent o\`u les $j$ sont des immersions ouvertes et $S_{\alpha} \rightarrow \mathcal{T}_{\alpha}$ une immersion.\\

On a un diagramme commutatif dont les carr\'es verticaux sont cart\'esiens, de m\^eme que le carr\'e horizontal  $\ \rondI$  en bas \`a droite\\

$$
\begin{array}{c}
\xymatrix{
& X^{(q/S_{\alpha})} \ar@{.>}[dd]^(.3){f^{(q)}_{\alpha}} |\hole \ar@{^{(}->}[rr]& & \tilde{X}'_{\alpha} \ar@{.>}[dd]^(.3){\tilde{f}'_{\alpha}} |\hole \ar[rr]& & \overline{\mathcal{Y}}'_{\alpha} \ar[dd]^{\overline{h}'_{\alpha}} \\
X^{(q/S_{\alpha})} \ar@{^{(}->}[rr] \ar[dd]_{f^{(q)}_{\alpha}} \ar@{=} [ur]^{\textrm{id}} & & \tilde{X}''_{\alpha} \ar[rr] \ar[dd]_(.3){\tilde{f}''_{\alpha}} \ar[ur]  & & \overline{\mathcal{Y}}''_{\alpha} \ar[dd]^(.3){\overline{h}''_{\alpha}} \ar[ur]& \\
& S_{\alpha} \ar@{^{(}.>}[rr]^(.7){j_{\tilde{S}'_{\alpha}}} |\hole & &\tilde{S}'_{\alpha}\ar@{^{(}.>}[rr]^(.4){i_{\tilde{S}'_{\alpha}}}  |\hole & &  \overline{\mathcal{T}}'_{\alpha}  \\
S_{\alpha} \ar@{^{(}->}[rr]_{j_{\tilde{S}''_{\alpha}}} \ar@{==}[ur]^{\textrm{id}} & & \tilde{S}'' _{\alpha} \ar@{^{(}->}[rr]_{i_{\tilde{S}''_{\alpha}}} \ar@{.>}[ur] \ar@{}[urrr] |{\rondI} & & \  \overline{\mathcal{T}}''_{\alpha} \ .\ar[ur]_{v_{\overline{\mathcal{T}}'_{\alpha} }} 
}
\end{array}
\leqno{(2.1.9)}
$$

\noindent D'apr\`es [Et 7, (3.4.4)] ou [Et 5, chap II, (3.4.4)], $v^{\ast}_{\overline{\mathcal{T}'}_{\alpha} K}$ induit un isomorphisme\\

\noindent (2.1.10)$ \qquad R^i f^{(q)}_{\alpha\  \textrm{rig}^{\ast}}(X^{(q/S_{\alpha})}_{\alpha} / \overline{\mathcal{T}}'_{\alpha}, \pi^{\ast}_{X_{\alpha}}\  j^{\ast}_{X_{\alpha}} E) \tilde{\rightarrow}  R^i f^{(q)}_{\alpha\  \textrm{rig}^{\ast}} (X^{(q/S_{\alpha})} / \overline{\mathcal{T}}''_{\alpha},  \pi^{\ast}_{X_{\alpha}}\  j^{\ast}_{X_{\alpha}} E) : $\\

\noindent en effet on peut appliquer [loc.cit.] car le morphisme propre $v_{\overline{\mathcal{T'}}_{\alpha}}$, \'etant \'etale au voisinage de $S_{\alpha}$, il induit [B 3, (1.3.5)] un isomorphisme entre un voisinage strict de $]S_{\alpha}[_{\overline{\mathcal{T}}''_{\alpha}}$ dans $\overline{\mathcal{T}}''_{\alpha K}$ et un voisinage strict de $]S_{\alpha}[_{\overline{\mathcal{T}}'_{\alpha}}$ dans $\overline{\mathcal{T}}'_{\alpha K}$.\\
Notons $\overline{T}''_{\alpha}$ la r\'eduction de $\overline{\mathcal{T}}''_{\alpha}\  \textrm{mod}\  \mathfrak{m}, v_{\overline{T}_{\alpha}}$ la r\'eduction de $v_{\overline{\mathcal{T}}_{\alpha}} \  \textrm{mod}\  \mathfrak{m}$ et $\tilde{S}'''_{\alpha }$ l'adh\'erence sch\'ematique de $S_{\alpha}$ dans ${\overline{T''}_{\alpha}}$ ; on a un diagramme commutatif o\`u les carr\'es $\  \rondI$ et $\ \rondII$  sont cart\'esiens\\

$$
\begin{array}{c}
\xymatrix{
& & \tilde{S}''_{\alpha} \ar@{^{(}->}[rrdd]^{i''} & & &&\\
& \tilde{S}'''_{\alpha} \ar@{^{(}->}[rd]^{i'''} \ar@{^{(}->}[ur]^{i'} & &&& &\\
& & v^{-1}_{\overline{T}_{\alpha}}(\tilde{S}_{\alpha}) \ar@{^{(}->}[rr] \ar[d]^{v_{\tilde{S}_{\alpha}}} \ar @{}[drr]|{\rondI}&& \overline{T}''_{\alpha} \ar@{^{(}->}[rr] \ar[d]^{v_{\overline{T}_{\alpha}}} \ar@{}[drr]|{\rondII}&& \overline{\mathcal{T}}''_{\alpha} \ar[d]^{v_{\overline{\mathcal{T}}_{\alpha}}} \\
S_{\alpha} \ar@{^{(}->}[rr]_{j_{\tilde{S}_{\alpha}}} \ar@{^{(}->}[uur]^{j_{\tilde{S}'''_{\alpha}}} & & \tilde{S}_{\alpha} \ar@{^{(}->}[rr]_{i_{\overline{T}_{\alpha}}}  && \overline{T}_{\alpha} \ar@{^{(}->}[rr]  && \overline{\mathcal{T}}_{\alpha}
}
\end{array}
\leqno{(2.1.11)}
$$

\noindent o\`u les $j$ (resp. les $i$) sont des immersions ouvertes (resp. ferm\'ees). Posons $v_{\alpha} = v_{\tilde{S}_{\alpha}} \circ i'''$.\\

\noindent Soit $\tilde{f}'''_{\alpha} : \tilde{X}'''_{\alpha} \rightarrow \tilde{S}'''_{\alpha}$ l'image inverse de $\tilde{f}''_{\alpha} : \tilde{X}''_{\alpha} \rightarrow \tilde{S}''_{\alpha}$ par $i' : \tilde{S}'''_{\alpha} \hookrightarrow \tilde{S}''_{\alpha}$. On a un diagramme commutatif \\

$$
\begin{array}{c}
\xymatrix{
X^{(q/S_{\alpha})} \ar@{^{(}->}[r] \ar[d]_{f^{(q)}_{\alpha}} & \tilde{X}'''_{\alpha} \ar@{^{(}->}[r] \ar[d]^{\tilde{f}'''_{\alpha}}  & \overline{\mathcal{Y}''}_{\alpha} \ar[d]^{\overline{h}''_{\alpha}} \\
S_{\alpha} \ar@{^{(}->}[r]^{j_{\tilde{S}'''_{\alpha}}} \ar@{=}[d]  & \tilde{S}'''_{\alpha} \ar[d]^{i_{\mathcal{T}}} \ar@{^{(}->}[r]^{i_{\tilde{S}'''_{\alpha}}} & \overline{\mathcal{T}''_{\alpha}} \ar[d]^{v_{\overline{\mathcal{T}''_{\alpha}}}}\\
S_{\alpha} \ar@{^{(}->}[r]^{j_{\tilde{S}_{\alpha}}} & \tilde{S}_{\alpha} \ar@{^{(}->}[r]^{i_{\tilde{S}_{\alpha}}} & \overline{\mathcal{T}_{\alpha}}
}
\end{array}
\leqno{(2.1.12)}
$$

\noindent o\`u les $j$ (resp. les $i$) sont des immersions ouvertes (resp. ferm\'ees). Comme $v_{\overline{\mathcal{T}}_{\alpha}}$ est \'etale au voisinage de $S_{\alpha}$ et que $v_{\alpha}$ est propre on d\'eduit de [B 3, th\'eo (1.3.5)] que $v_{\overline{\mathcal{T}}_{\alpha K}}$ induit un isomorphisme entre un voisinage strict de $]S_{\alpha}[_{\overline{\mathcal{T}}''_{\alpha}}$ dans $]\tilde{S}'''_{\alpha}[_{\overline{\mathcal{T}}''_{\alpha}}$ et un voisinage strict de $]S_{\alpha}[_{\overline{\mathcal{T}}_{\alpha}}$ dans  $]\tilde{S}_{\alpha}[_{\overline{\mathcal{T}}_{\alpha}}$. Par suite [Et 7, th\'eo (3.4.4)] $v_{\overline{\mathcal{T}}_{\alpha K}}$ induit un isomorphisme\\

\noindent (2.1.13) $ R^i f^{(q)}_{\alpha\  \textrm{rig}^{\ast}}(X^{(q/S_{\alpha})}_{\alpha} / \overline{\mathcal{T}''}_{\alpha}, \pi^{\ast}_{X_{\alpha}}\  j^{\ast}_{X_{\alpha}} E) \tilde{\rightarrow}  R^i f^{(q)}_{\alpha\  \textrm{rig}^{\ast}} (X^{(q/S_{\alpha})} / \overline{\mathcal{T}}_{\alpha},  \pi^{\ast}_{X_{\alpha}}\  j^{\ast}_{X_{\alpha}} E).$\\

Par composition des isomorphismes (2.1.8), (2.1.10) et (2.1.13) on obtient un isomorphisme induit par $\pi_{X_{\alpha}} : X^{(q/S_{\alpha})}_{\alpha} \rightarrow X_{\alpha}$

$$ 
\hspace{-1cm}\begin{array}{c}
\xymatrix{
F^{\ast}_{S_{\alpha}}\, R^i\, f_{\alpha \textrm{rig}^{\ast}}(X_{\alpha}/ \overline{\mathcal{T}_{\alpha}}; j^{\ast}_{X_{\alpha}}\, E) \ar[r]^{\sim} \ar[d]^{\simeq} &  R^i\, f^{(q)}_{\alpha \textrm{rig}^{\ast}}(X^{(q/S_{\alpha}}/ \overline{\mathcal{T}}_{\alpha}; \pi^{\ast}_{X_{\alpha}} j^{\ast}_{X_{\alpha}} E)\\
F^{\ast}_{\overline{\mathcal{T}}'_{\alpha K}} R^i f_{\alpha \textrm{rig}^{\ast}}(X_{\alpha}/\overline{\mathcal{T}_{\alpha}}; j^{\ast}_{X_{\alpha}} E) . & 
}
\end{array}
\leqno{(2.1.14)}
$$

Si l'on prend l'image inverse de cet isomorphisme par $j_{\overline{\mathcal{T}}'_{\alpha}} : \mathcal{T}_{\alpha} \rightarrow \overline{\mathcal{T}}'_{\alpha}$, on voit ais\'ement en suivant les diagrammes pr\'ec\'edents que l'on obtient l'isomorphisme en cohomologie convergente, induit par $\pi_{X_{\alpha}}$ :\\

\noindent (2.1.15) $\quad F^{\ast}_{S_\alpha} R^i f_{\alpha\  \textrm{conv}^{\ast}} (X_{\alpha} / \mathcal{T}_{\alpha}, j^{\ast}_{X_{\alpha}} \mathcal{E})\  \tilde{\rightarrow}\  R^{i} f^{(q)}_{\alpha\  \textrm{conv}^{\ast}} (X^{(q/S_{\alpha}}) / \mathcal{T}_{\alpha},\pi^{\ast}_{X_{\alpha}}\  j^{\ast}_{X_{\alpha}} \mathcal{E})$\\

\noindent o\`u l'on a, comme pour (2.1.5), un isomorphisme\\

\noindent (2.1.16) $\qquad \qquad j^{\ast}_{S_{\alpha}} R^i f_{\textrm{conv}^{\ast}} (X / \mathcal{S}, \mathcal{E})\  \tilde{\rightarrow}\  R^i f_{\alpha\ \textrm{conv}^{\ast}} (X_{\alpha} / \mathcal{T}_{\alpha}, j^{\ast}_{X_{\alpha}} \mathcal{E}).$\\

Si $F_{\tilde{S}_{\alpha}}$ d\'esigne le Frobenius (puissance $q$) de $\tilde{S}_{\alpha}$, on notera $\tilde{X}_{\alpha}^{(q/\tilde{S}_{\alpha})}$ le produit fibr\'e  d\'efini par le diagramme \`a carr\'e cart\'esien

$$
\xymatrix{
\tilde{X}_{\alpha} \ar[r]^(.4){F_{\tilde{X}_{\alpha}/\tilde{S}_{\alpha}}} \ar[rd]_{\tilde{f}_{\alpha}} \ar@/^4pc/[rr]^{F_{\tilde{X}_{\alpha}}} & \tilde{X}_{\alpha}^{(q/\tilde{S}_{\alpha})} \ar[r]^{\pi_{\tilde{X}_{\alpha}}} \ar[d]^{\tilde{f}^{(q)}_{\alpha}} & \tilde{X}_{\alpha}  \ar[d]^{\tilde{f}_{\alpha}} \\
& \tilde{S}_{\alpha} \ar[r]_(.48){F_{\tilde{S}_{\alpha}}}  \ar@{^{(}->}[d]_{i_{\tilde{S}_{\alpha}}} & \tilde{S}_{\alpha}\\
& \overline{\mathcal{T}_{\alpha}} &.
}
$$

\noindent Le calcul de $R^i f^{(q)}_{\alpha\  \textrm{rig}^{\ast}}(X^{(q/S_{\alpha})} / \overline{\mathcal{T}}_{\alpha}, \pi^{\ast}_{X_{\alpha}}\  j^{\ast}_{X_{\alpha}} E)$ \'etant ind\'ependant de la compactification de $X^{(q/S_{\alpha}})$ choisie [B 5, (3.1.11), (3.1.12), (3.2.3)], nous choisirons dor\'enavant  $\tilde{X}_{\alpha}^{(q/ \tilde{S}_{\alpha})}$ [C-T, \S\ 10] comme compactification de $X^{(q/s_{\alpha})}$ au lieu de $\tilde{X}'''_{\alpha}$. Consid\'erons le diagramme commutatif

$$
\xymatrix{
X_{\alpha} \ar[rr]^{F_{X_{\alpha}/S_{\alpha}}} \ar@{^{(}->}[d]_{j_{\tilde{X}_{\alpha}}} && X_{\alpha}^{(q/S_{\alpha})} \ar[d]  \\
\tilde{X}_{\alpha} \ar[rr]^{F_{\tilde{X}_{\alpha}/\tilde{S}_{\alpha}}} \ar[d]_{i_{\tilde{S}_{\alpha}} \circ \tilde{f}_{\alpha}} && \tilde{X}^{(q/\tilde{S}_{\alpha})}_{\alpha}    \ar[d]^{i_{\tilde{S}_{\alpha}} \circ \tilde{f}^{(q)}_{\alpha}}   \\
\overline{\mathcal{T}}_{\alpha} \ar@{=}[rr] && \overline{\mathcal{T}}_{\alpha}   ;  \\
} 
$$

\noindent $F_{X_{\alpha}/S_{\alpha}}$ d\'efinit par fonctorialit\'e [B 5, (3.1.11) (ii), (3.1.12) (i)] ou [C-T, (10.5.2)]  un morphisme\\

$$
\begin{array}{c}
\xymatrix{
F^{\ast}_{X_{\alpha}/S_{\alpha}} :  R^i f^{(q)}_{\alpha\  \textrm{rig}^{\ast}} (X^{(q/S_{\alpha})} / \overline{\mathcal{T}}_{\alpha}; \pi^{\ast}_{X_{\alpha}}\  j^{\ast}_{X_{\alpha}} E) \ar[r] & R^i f_{\alpha\  \textrm{rig}^{\ast}}(X_{\alpha} / \overline{\mathcal{T}}_{\alpha}; F^{\ast}_{X_{\alpha}/S_{\alpha}}  \pi^{\ast}_{X_{\alpha}}  j^{\ast}_{X_{\alpha}} E) \ar[d]^{\simeq}\\
  R^i f_{\alpha\  \textrm{rig}^{\ast}} (X_{\alpha}/ \overline{\mathcal{T}}_{\alpha};  j^{\ast}_{X_{\alpha}} F^{\ast}_{X} E)\ar[r]^{ \simeq} &R^i f_{\alpha\  \textrm{rig}^{\ast}} (X_{\alpha}/  \overline{\mathcal{T}}_{\alpha}; F^{\ast}_{X_{\alpha}}\  j^{\ast}_{X_{\alpha}} E) .
}
\end{array}
\leqno{(2.1.17)}
$$

Par image inverse par $j_{\overline{\mathcal{T}}_{\alpha}} : \mathcal{T}_{\alpha} \hookrightarrow \overline{\mathcal{T}}_{\alpha}$ il fournit le morphisme en cohomologie convergente, induit par $F_{X_{\alpha}/S_{\alpha}}$\\

\noindent (2.1.18) $F^{\ast}_{X_{\alpha}/S_{\alpha}} : R^i f^{(q)}_{\alpha\  \textrm{conv}^{\ast}} (X^{(q/S_{\alpha})} / \mathcal{T}_{\alpha}, \pi^{\ast}_{X_{\alpha}}\  j^{\ast}_{X_{\alpha}} \mathcal{E}) \rightarrow R^i f_{\alpha\  \textrm{conv}^{\ast}}(X_{\alpha} / \mathcal{T}_{\alpha}, F^{\ast}_{X_{\alpha}}\  j^{\ast}_{X_{\alpha}} \mathcal{E}).$ \\

Enfin, l'isomorphisme de Frobenius de $E, \phi_{E} : F^{\ast}_{X} \tilde{\rightarrow} E$, fournit par fonctorialit\'e un isomorphisme\\

\noindent (2.1.19) $\qquad  R^i f_{\alpha\  \textrm{rig}^{\ast}} (X_{\alpha} / \overline{\mathcal{T}}_{\alpha}, j^{\ast}_{X_{\alpha}}  F^{\ast}_{X} E)\  \tilde{\rightarrow}\  R^{i} f_{\alpha\  \textrm{rig}^{\ast}} (X_{\alpha} / \overline{\mathcal{T}}_{\alpha},    j^{\ast}_{X_{\alpha}} \mathcal{E}),$\\

\noindent dont l'image inverse par $j_{\overline{\mathcal{T}}_{\alpha}} : \mathcal{T}_{\alpha} \rightarrow \overline{\mathcal{T}}_{\alpha}$ est l'isomorphisme induit en cohomologie convergente par $\phi_{\mathcal{E}} : F^{\ast}_{X}\  \mathcal{E}\  \tilde{\rightarrow}\  \mathcal{E}$,\\

\noindent (2.1.20) $\qquad R^i f_{\alpha\  \textrm{conv}^{\ast}} (X_{\alpha} / \mathcal{T}_{\alpha}, j^{\ast}_{X_{\alpha}}  F^{\ast}_{X} \mathcal{E})\  \tilde{\rightarrow}\  R^{i} f_{\alpha\  \textrm{conv}^{\ast}} (X_{\alpha} / \mathcal{T}_{\alpha},    j^{\ast}_{X_{\alpha}} \mathcal{E}).$\\

En composant les morphismes (2.1.14), (2.1.17) et (2.1.19) [resp. (2.1.15), (2.1.18) et (2.1.20)] on obtient le morphisme de Frobenius souhait\'e\\

\noindent (2.1.21) $\qquad  \phi_{E_{\alpha_{i}}} : F^{\ast}_{S_{\alpha}} R^i f_{\alpha\  \textrm{rig}^{\ast}} (X_{\alpha} / \overline{\mathcal{T}}_{\alpha},   j^{\ast}_{X_{\alpha}} E)\  \rightarrow\  R^i f_{\alpha\  \textrm{rig}^{\ast}}(X_{\alpha} / \overline{\mathcal{T}}_{\alpha}, \   j^{\ast}_{X_{\alpha}} E)$\\

\noindent [resp.\\

\noindent (2.1.22) $\qquad  \phi_{\mathcal{E}_{\alpha_{i}}} : F^{\ast}_{S_{\alpha}} R^i f_{\alpha\  \textrm{conv}^{\ast}} (X_{\alpha} / \mathcal{T}_{\alpha},   j^{\ast}_{X_{\alpha}} \mathcal{E})\  \rightarrow\  R^i f_{\alpha\  \textrm{conv}^{\ast}}(X_{\alpha} / \mathcal{T}_{\alpha}, \   j^{\ast}_{X_{\alpha}} \mathcal{E})$\\

\noindent qui est l'image inverse de $\phi_{E_{\alpha_{i}}}$ par $j_{\overline{\mathcal{T}}_{\alpha}}].$\\

D'apr\`es [Et 8, th\'eo (3.3.1)] ou [Et 5, chap III, th\'eo (3.3.1)] $\phi_{\mathcal{E}_{\alpha_{i}}}$ est un isomorphisme : nous allons en d\'eduire que $\phi_{E_{\alpha_{i}}}$ est un isomorphisme, ce qui ach\`evera la preuve du th\'eor\`eme.\\

On a un diagramme commutatif\\

$$
\begin{array}{c}
\xymatrix{
& \tilde{S}_{\alpha} \ar[d]^{i_{\overline{T}_{\alpha}}} \ar@{^{(}->}[rd]^{i_{\tilde{S}_{\alpha}}}  \ar@/_1pc/[dd]_{v} &\\
S_{\alpha} \ar@{^{(}->}[ru]^{j_{\tilde{S}_{\alpha}}}    \ar@{^{(}->}[rd]_{j_{\overline{S}_{\alpha}}} & \overline{T}_{\alpha} \ar[d]^{\overline{u}'_{\alpha}} \ar@{^{(}->}[r]_{i_{\overline{T}_{\alpha}}}   & \overline{\mathcal{T}}_{\alpha} \ar[d]^{\overline{u}_{\alpha}}\\
& \overline{S}_{\alpha} \ar@{^{(}->}[r]_{i_{\overline{S}_{\alpha}}} & \overline{\mathcal{S}}_{\alpha}
}
\end{array}
\leqno{(2.1.23)}
$$

\noindent o\`u $j_{\overline{S}_{\alpha}}$ (resp. $\overline{u}'_{\alpha}$) est la r\'eduction mod $\mathfrak{m}$ de $j_{\overline{\mathcal{S}}_{\alpha}} : \mathcal{S}_{\alpha} \hookrightarrow \overline{\mathcal{S}}_{\alpha}$ (resp. de $\overline{u_{\alpha}} : \overline{\mathcal{T}_{\alpha}} \rightarrow \overline{\mathcal{S}_{\alpha}})$ : les $j$ (resp. les $i$) sont des immersions ouvertes (resp. ferm\'ees).\\

De m\^eme le triangle commutatif\\

$$
\xymatrix{
& \mathcal{T}_{\alpha} = \mathcal{S}_{\alpha} \times_{\mathcal{V}} \mathcal{S} \ar[dd]^{u_{\alpha}} \\
\mathcal{S}_{\alpha} \ar@{^{(}->} [ur]^{\Delta} \ar@{=}[rd]_{\textrm{id}}\\
& \mathcal{S}_{\alpha}
}
$$

\noindent o\`u $\Delta$ est le morphisme diagonal et $u_{\alpha}$ la projection ($u_{\alpha}$ est lisse), fournit un triangle commutatif\\

$$
\begin{array}{c}
\xymatrix{
&   \mathcal{T}_{\alpha} \ar[d]^{u_{\alpha}} \\
S_{\alpha} \ar@{^{(}->}[ur]^{i'} \ar@{^{(}->}[r]_{i} & \mathcal{S}_{\alpha}\\
}
\end{array}
\leqno{(2.1.24)}
$$

\noindent o\`u $i$ et $i'$ sont des immersions ferm\'ees.\\
D'apr\`es [B 3, (2.3.1)] le foncteur $u^{\ast}_{\alpha_{K}}$ induit une auto-\'equivalence de la cat\'egorie $F^a\mbox{-}\textrm{Isoc}(S_{\alpha}/K)$ : en composant un foncteur quasi-inverse canonique \`a $u^{\ast}_{K}$ (cf. [B 5, (3.1.10)]) avec l'\'equivalence de cat\'egories de [Et 8, cor (1.2.3)] ou [Et 5, chap III, cor (1.2.3)] on constate que la donn\'ee de l'isomorphisme $\phi_{\mathcal{E}_{\alpha_{i}}}$ correspond \`a la donn\'ee d'un isomorphisme 

$$
\phi_{\mathcal{M}} : \mathcal{M}^{\sigma} \tilde{\longrightarrow}\  \mathcal{M}
$$

\noindent de $\hat{A}_{K}$-modules projectifs de type fini commutant aux connexions.\\

De m\^eme, d'apr\`es [B 3, (2.3.5)] le morphisme propre $v$ de (2.1.23) induit l'auto-\'equivalence $v^{\ast} = \overline{u}^{\ast}_{\alpha_{k}}$ de Isoc$^{\dag}(S_{\alpha}/K)$ : en composant un foncteur quasi-inverse \`a $v^{\ast}$ avec l'\'equivalence de cat\'egories de [B 3, (2;5;2) (ii)], on constate que la donn\'ee du morphisme $\phi_{E_{\alpha_{i}}}$ correspond \`a la donn\'ee d'un morphisme

$$
\phi_M : M^{\sigma} \longrightarrow M
$$

\noindent de $A^{\dag}_{K}$-modules projectifs de type fini commutant aux  connexions.\\

Ce qui pr\'ec\`ede peut \^etre formalis\'e par un diagramme commutatif de foncteurs entre cat\'egories\\

$$
\begin{array}{c}
\xymatrix{
\textrm{Conn}^{\dag} (A^{\dag}_{K}) \ar[d]^{\mathcal{G}} && \textrm{Isoc}^{\dag} (S_{\alpha}/K) \ar[ll]_{\Gamma(\overline{\mathcal{S}}_{\alpha_{K},-})} ^{\simeq} \ar[rr]^{\overline{u}^{\ast}_{\alpha_K}}_{\simeq} \ar[d]^{j^{\ast}_{\overline{\mathcal{S}}_{\alpha K}}} && \textrm{Isoc}^{\dag}(S_{\alpha}/K) \ar[d]^{j^{\ast}_{\overline{\mathcal{T}}_{\alpha K}}} \\
\textrm{Conn}^{\wedge} (\hat{A}_{K}) && \textrm{Isoc}(S_{\alpha}/K) \ar[ll]_{\Gamma(\mathcal{S}_{\alpha_{K},-})} ^{\simeq}   \ar[rr]^{u^{\ast}_{\alpha_K}}_{\simeq}  && \textrm{Isoc}(S_{\alpha}/K)
}
\end{array}
\leqno{(2.1.25)}
$$

\noindent o\`u les fl\`eches verticales sont les foncteurs d'oubli et les fl\`eches horizontales des \'equivalences de cat\'egories : pour les notations et r\'esultats cf. [Et 8, (1.1) et prop (1.2.1)] ou [Et 5, chap III, (1.1) et prop. (1.2.1)]. Ainsi on a un carr\'e commutatif 

$$
\xymatrix{
\mathcal{G}(\phi_{M}) : \mathcal{G}(M)^{\sigma} \ar[r] \ar[d]^{\simeq} & \mathcal{G}(M) \ar[d]^{\simeq}\\
\phi_{\mathcal{M}} : \mathcal{M}^{\sigma} \ar[r]^{\sim} & \mathcal{M}\ ,
 }
$$

\noindent  donc $\mathcal{G}(\phi_{M})$ est un isomorphisme : par fid\`ele platitude de $\hat{A}_{K}$ sur $A^{\dag}_{K}$ on en d\'eduit que $\phi_{M}$ est un isomorphisme. Par suite $\phi_{E_{\alpha_{i}}}$ est un isomorphisme. $\square$

\newpage
 \noindent \textbf{Remarque (2.2)}. 
\begin{itemize}
\item[(i)] En fait on a prouv\'e, plus pr\'ecis\'ement, que le morphisme (2.1.17) induit par $F_{X_{\alpha}/S_{\alpha}}$ est un isomorphisme.
\item[(ii)] Pour construire le morphisme de Frobenius $\phi_{E_{i}}$ nous n'avons pas suppos\'e l'existence d'un rel\`evement \`a $\overline{\mathcal{S}}$ du Frobenius de $\overline{S}$, contrairement \`a [C-T, 12.2]. 
\end{itemize}

\vskip 3mm
\section*{3. Cas propre et lisse}

\noindent \textbf{Th\'eor\`eme (3.1)}. \textit{Soient $S$ un $k$-sch\'ema lisse et s\'epar\'e et $f : X \rightarrow S$ un $k$-morphisme projectif et lisse satisfaisant aux propri\'et\'es de [Et 7, (3.4.8.2), ou (3.4.8.6), ou  (3.4.9)] ( cf [Et 5, chap II, (3.4.8.2), ou (3.4.8.6), ou  (3.4.9)]). Alors} \\

\noindent (3.1.1) \textit{Pour tout entier $i \geqslant 0$, on a un diagramme commutatif de foncteurs naturels induits par $f$ et d\'efinis dans [Et 7, (3.4.8.5)] ou [Et 5, chap II, (3.4.8.5)]}
$$
\xymatrix{
 F^{a}\mbox{-}Isoc^{\dag}(X/ K) \ar[rr]^{R^{i}f_{rig\ast}}\ar[d]&&F^{a}\mbox{-}Isoc^{\dag}(S/K)\ar [d]\\
F^{a}\mbox{-}Isoc (X/ K)\ar[rr]^{R^{i}f_{conv\ast}}&&F^{a}\mbox{-}Isoc(S/ K)
}
$$
\textit{o\`u les fl\`eches verticales sont les foncteurs d'oubli.}\\

\noindent (3.1.2) \textit{Le foncteur $R^i f_{\textrm{rig}^{\ast}}$ pr\'ec\'edent est compatible aux changements de base entre $k$-sch\'emas lisses et s\'epar\'es (en particulier il commute aux passages aux fibres en les points ferm\'es de $S$), c'est-\`a-dire : pour tout carr\'e cart\'esien}

$$
\xymatrix{
X' \ar[r]^{g'} \ar[d]_{f'} & X \ar[d]^f\\
S' \ar[r]_{g} & S
}
$$

\noindent \textit{o\`u $S'$ est un $k$-sch\'ema lisse et s\'epar\'e et $E  \in F^a\mbox{-}\textrm{Isoc}^{\dag}(X/K)$ on a un isomorphisme de changement de base}

$$g^{\ast} R^i f_{\textrm{rig}^{\ast}}(E) \tilde{\longrightarrow} R^i f'_{\textrm{rig}^{\ast}}(g'^{\ast} (E)) $$

\noindent \textit{compatible aux connexions et aux Frobenius.}

\vskip 3mm
\noindent \textit{D\'emonstration}. \\
\noindent \textit{Pour (3.1.1)}. Vu la d\'efinition locale sur $S$ de $R^i f_{\textrm{rig}^{\ast}}$ [cf [Et 7, (3.4.8.5)]) on peut supposer $S$ affine lisse et connexe et se ramener au cas de [Et 7, (3.4.8.2)]: alors $f$ est relevable comme dans le th\'eor\`eme (2.1) ci-dessus qu'il suffit d'appliquer, d'o\`u la conclusion.\\

\noindent \textit{Pour (3.1.2)}. Soient $V = Spec\ A_{0} \displaystyle \mathop{\hookrightarrow}^{j_{U_{0}}}S$ un ouvert affine de $S$ et $Y = Spec\ B_{0} \displaystyle \mathop{\hookrightarrow}_{j'_{Y_{0}}}S'$ un ouvert affine de $V' = S' \times_{S} V$ ; on note $\psi_{VY} : Y \rightarrow V$, le morphisme induit par $g$. Soient $A = \mathcal{V}[t_{1},...,t_{n}]/(f_{1},...f_{r})$ une $\mathcal{V}$-alg\`ebre lisse relevant $A_{0}$, $U = Spec\ A$ et $\overline{U}$ la fermeture projective de $U$ dans $\mathbb{P}^n_{\mathcal{V}}$, $\mathcal{S}$ et $\overline{\mathcal{S}}$ leurs compl\'et\'es formels respectifs, $j_{\overline{\mathcal{S}}}$ l'immersion ouverte $\mathcal{S} \hookrightarrow \overline{\mathcal{S}}$, et posons $X_{V} = X \times_{S} V$, et $\overline{V} = \overline{U}\  \textrm{mod} . \pi$. Si $V$ (resp. $Y$) parcourt un recouvrement ouvert affine de $S$ (resp. de $V'$) alors $Y$ parcourt un recouvrement ouvert affine de $S'$ ; or la donn\'ee de $R^i f_{\textrm{rig}^{\ast}}(E)$ (resp. de $g^{\ast} R^i f_{\textrm{rig}^{\ast}}(E))$ \'equivaut \`a la donn\'ee des $j^{\ast}_{V} R^i f_{\textrm{rig}^{\ast}}(E)$ (resp. des $j_{Y}^{\prime\ast} g^{\ast} R^i f_{\textrm{rig}^{\ast}}(E))$, donc la donn\'ee de $g^{\ast} R^i f_{\textrm{rig}^{\ast}}(E)$ \'equivaut \`a celle des

$$
\psi^{\ast}_{VY}\  j^{\ast}_{V}\  R^i f_{\textrm{rig}^{\ast}}(E) = \psi^{\ast}_{VY}\  R^i f_{\textrm{rig}^{\ast}}(X_{V} / \overline{\mathcal{S}}, E_{X_{V}}).
$$
Puisque $\psi  := \psi_{VY}$ est de type fini on peut choisir une pr\'esentation $B_{0} = A_{0}[x_{1},...,x_{d}]/(g_{1},...,g_{s})$ de $B_{0}$ sur $A_{0}$ : notons $\mathcal{Y}$ (resp. $\overline{\mathcal{Y}}$) le compl\'et\'e formel de $\mathbb{A}^d_{U}$ (resp. de $\mathbb{P}^d_{\overline{U}})$, $\overline{Y}$ l'adh\'erence sch\'ematique de $Y$ dans $\mathbb{P}^d_{\overline{U}}$, $\overline{\psi} : \overline{Y} \rightarrow \overline{V}$ le morphisme canonique et $j_{\overline{\mathcal{Y}}} : \mathcal{Y} \hookrightarrow \overline{\mathcal{Y}}$, $j_{\overline{Y}} : Y \hookrightarrow \overline{Y}$ les immersions ouvertes ; $\theta : \mathcal{Y} \rightarrow \mathcal{S}$, $\overline{\theta} : \overline{\mathcal{Y}} \rightarrow \overline{\mathcal{S}}$ d\'esignerons les projections canoniques. D'o\`u un diagramme commutatif

$$
\xymatrix{
& \overline{Y} \ar@{^{(}->}[rr] \ar@{.>}[dd]^(.4){\overline{\psi}} |\hole && \overline{\mathcal{Y}} \ar[dd]^{\overline{\theta}} \\
Y \ar@{^{(}->}[rr] \ar[dd]_{\psi} \ar@{^{(}->}[ur] && \mathcal{Y} \ar[dd] \ar@{^{(}->}[ur]_{j_{\overline{\mathcal{Y}}}}\\
&\overline{V} \ar@{.>}[rr]^(.40){\theta}  && \overline{\mathcal{S}}  \\
V \ar@{^{(}->}[rr] \ar@{^{(}.>}[ur] &&\  \mathcal{S}\  .\ar@{^{(}->}[ur]_{j_{\overline{\mathcal{S}}}} \\
}
$$
\\

Ainsi [B 3, (2.3.2) (iv)]

$$
\psi^{\ast} R^i f_{\textrm{rig}^{\ast}}(X_{V}/ \overline{\mathcal{S}}, E_{X_{V}}) = \overline{\theta}^{\ast} R^i f_{\textrm{rig}^{\ast}} (X_{V}/ \overline{\mathcal{S}}, E_{X_{V}})
$$

\noindent et \quad $j^{\ast}_{\overline{\mathcal{Y}}}\  \overline{\theta}^{\ast} R^i f_{\textrm{rig}^{\ast}}(X_{V} / \overline{\mathcal{S}}, E_{X_{V}})$\\

$\quad =  \theta^{\ast} j^{\ast}_{\overline{\mathcal{S}}}\ R^i f_{\textrm{rig}^{\ast}}(X_{V} / \overline{\mathcal{S}}, E_{X_{V}})$\\

$\quad = \theta^{\ast} R^i f_{\textrm{conv}^{\ast}}(X_{V} / S, \hat{E}_{X_{V}})\ [\mbox{Et 7},(3.4.4)] \ \mbox{ou} \  [\mbox{Et 5, chap II}, (3.4.4)]$\\

$\quad = \psi^{\ast} R^i f_{\textrm{conv}^{\ast}}(X_{V} / S, \hat{E}_{X_{V}}) = \psi^{\ast} (R^i f_{\textrm{conv}^{\ast}} (\hat{E})_{\vert_{V}})$\\

$\qquad \qquad \qquad \qquad \qquad \qquad \quad = (g^{\ast} (R^i f_{\textrm{conv}^{\ast}} (\hat{E}))_{\vert_{Y}}$\\

\noindent o\`u $\hat{E}$ est l'isocristal convergent associ\'e \`a $E$.\\

De m\^eme on a\\

$j^{\ast}_{\overline{\mathcal{Y}}}\ j'^{\ast}_{Y}\  R^i f'_{\textrm{rig}^{\ast}} (g'^{\ast} (E)) \simeq R^i f'_{\textrm{conv}^{\ast}} (X'_{Y}/\mathcal{Y}, g'^{\ast}(\hat{E})_{X'_{Y}})$\\

$\qquad \qquad \qquad \qquad \qquad \simeq (R^i f'_{\textrm{conv}^{\ast}} (g'^{\ast}(\hat{E})))_{\vert_{Y}}$\ . \\

Or le th\'eor\`eme [Et 8, (3.2.1)] ou [Et 5, chap III, (3.2.1)] fournit un isomorphisme de changement de base

$$
g^{\ast} R^i f_{\textrm{conv}^{\ast}}(\hat{E})\  \tilde{\rightarrow}\ R^i f'_{\textrm{conv}^{\ast}} (g'^{\ast}(\hat{E}))\  ; 
$$

\noindent donc par le th\'eor\`eme de pleine fid\'elit\'e pour les $F$-isocristaux de Kedlaya [Ked 2, theo 1.1] on en d\'eduit un isomorphisme

$$
g^{\ast} R^i f_{\textrm{rig}^{\ast}}(E)\  \tilde{\rightarrow}\  R^i f'_{\textrm{rig}^{\ast}} g'^{\ast}(E)\ 
$$

\noindent compatible aux connexions et aux Frobenius. $\square$

\vskip 3mm
\noindent \textbf{Th\'eor\`eme (3.2)}. \textit{Soient $S = \displaystyle \mathop{\coprod}_{\alpha = 1}^n S_{\alpha}$ un $k$-sch\'ema lisse et s\'epar\'e, d\'ecompos\'e en la somme de ses composantes connexes, et $f : X \rightarrow S$ un $k$-morphisme propre et lisse v\'erifiant la propri\'et\'e $\mathcal{P}$ suivante :}\\

 \noindent $\mathcal{P}$
 $\left\lbrace
 \begin{array}{l}
 Il\  existe\  un\  ouvert\  dense\  U \subset X quasi\mbox{-}projectif sur\  S\ 
tel\  que\\ 
pour\  tout\  \alpha \in \mathbb{[[}1,n \mathbb{]]}\   on\  ait\  S_{\alpha} \backslash f(X \backslash U) \neq \phi . 
\end{array}
 \right .$\\

 \noindent \textit{Alors}\\
\begin{itemize} 
\item[{(3.2.1)}] \textit{La propri\'et\'e ($\mathcal{P})$ \'equivaut \`a dire que $f$ est g\'en\'eriquement projective, i.e. qu'il existe un ouvert dense $V \subset S$ tel que l'application $f_{V} : X_{V} = X_{X_{S}} V \rightarrow V$ induite par $f$ soit projective et lisse.\\
}

\item[{(3.2.2)}]\textit{Supposons $k$ parfait, $e\leqslant p-1$ et que le $f_{V}$ de (3.2.1) satisfait aux hypoth\`eses de [Et 7, (3.4.8.2), ou (3.4.8.6), ou (3.4.9)]}. $ Si\  E \in F^a\mbox{-}\textrm{Isoc}^{\dag}(X/K)_{\textrm{plat}} \ a\ pour\ image\ \mathcal{E} \in F^a\mbox{-}\textrm{Isoc}(X/K), alors :$

	\begin{itemize}
	\item[{(i)}] $\mathcal{E}^i = R^i   f_{\textrm{conv}^{\ast}}(\mathcal{E}) \in F^a\mbox{-}\textrm{Isoc}(S/K)$,

	\item[{(ii)}] \textit{Il existe $E^i \in F^a\mbox{-}\textrm{Isoc}^{\dag}(S/K)$, unique \`a isomorphisme pr\`es, tel que $\mathcal{E}^i$ soit l'image de $E^i$ par le foncteur d'oubli}

	$$F^a\mbox{-}\textrm{Isoc}^{\dag}(S/K) \longrightarrow F^a\mbox{-}\textrm{Isoc}(S/K) $$
	$$E^i \longmapsto \widehat{E^i} = \mathcal{E}^i . $$
	\end{itemize}
\end{itemize}

\vskip 3mm
\noindent \textit{D\'emonstration}. \\
\noindent \textit{Prouvons (3.2.1)}. Si $f$ est g\'en\'eriquement projective et lisse on prouve facilement qu'elle v\'erifie $(\mathcal{P})$.\\

R\'eciproquement supposons que $f$ v\'erifie $(\mathcal{P})$.\\

Par le lemme de Chow pr\'ecis de Gruson-Raynaud [R-G, I, cor 5.7.14] il existe un \'eclatement $U$-admissible $g : X' \rightarrow X$, avec $X'$ quasi-projectif sur $S$ : en particulier $g$ induit un isomorphisme

$$
g_{U} : U' = g^{-1}(U)\   \tilde{\longrightarrow}\  U\ .
$$

\noindent De plus, comme $f$ et $g$ sont propres, le morphisme compos\'e $f \circ g : X' \rightarrow X$ est projectif [EGA II, (5.5.3) (ii)]. L'image du ferm\'e $Z := X \backslash U$ par le morphisme propre $f$ est un ferm\'e $f(Z)$ de $S$ et l'ouvert $V = S \backslash f(Z) = \displaystyle \mathop{\coprod}_{\alpha} (S_{\alpha} \backslash f(X-U))$ est non vide par hypoth\`ese : comme $S_{\alpha}$ est connexe et int\`egre, l'ouvert non vide $V_{\alpha} := S_{\alpha} \backslash f(X  \backslash U)$ de $S_{\alpha}$ est dense, donc l'immersion ouverte $j : V \hookrightarrow S$ est dominante. D'autre part l'ouvert $X_{V} = X \times_{S} V$ de $X$ ne rencontre pas $f^{-1}(f(Z))$, donc $X_{V}$ est un ouvert de $U$ : par suite l'isomorphisme $g_{U}$ induit un isomorphisme 

$$
g_{V} : X'_{U} = g^{-1}(X_{V}) \tilde{\longrightarrow} X_{V} .
$$

\noindent Notons $f_{V} : X_{V} \rightarrow V$ la restriction de $f$ ; le morphisme compos\'e $f_{V} \circ g_{V}$, restriction du morphisme projectif $f \circ g$, est lui aussi projectif : ainsi $f_{V}$ est projectif, d'o\`u (3.2.1).\\

\noindent \textit{Prouvons le (3.2.2)}.  Le (i) est mis pour m\'emoire, car prouv\'e en [Et 8, (3.3.1)] ou [Et 5, chap III, (3.3.1)]. Pour le (ii) consid\'erons le carr\'e cart\'esien

$$
\xymatrix{
X_{V} \ar@{^{(}->}[r]^{j'} \ar[d]_{f_{V}} & X \ar[d]^{f}\\
V \ar@{^{(}->}[r]^{j} & S   ; 
}
$$

\noindent on a un isomorphisme de changement de base en cohomologie convergente [Et 8, (3.3.1)]\\

\noindent (3.2.2.1) $\qquad \qquad  j^{\ast} R^i   f_{\textrm{conv}^{\ast}}(\mathcal{E}) \tilde{\longrightarrow} R^i   f_{V{\textrm{conv}}^{\ast}}(j'^{\ast} (\mathcal{E})) =: \mathcal{E}^i_{V}, $\\

\noindent o\`u $\mathcal{E}$ d\'esigne l'image de $E$ par le foncteur d'oubli $F^a\mbox{-}\textrm{Isoc}^{\dag}(X/K) \rightarrow F^a\mbox{-}\textrm{Isoc}(X/K)$, $E \mapsto \hat{E} = \mathcal{E}$.\\
Pour la suite de la d\'emonstration on peut supposer $S$ connexe, int\`egre : quitte \`a restreindre $V$ on peut supposer $V$ affine, lisse et connexe, $V = Spec\ A_{0}$. On utilise alors les notations introduites dans la d\'emonstration du th\'eor\`eme (3.1) :  $j_{\overline{\mathcal{S}}} : \mathcal{S} = Spf\ \hat{A} \hookrightarrow \overline{\mathcal{S}}$.
De plus $f_{V}$ se rel\`eve en un morphisme projectif et lisse  $h : \mathcal{X} \rightarrow \mathcal{S}$ s'ins\'erant dans un carr\'e cart\'esien de $\mathcal{V}$-sch\'emas formels
 
$$
\xymatrix{
\mathcal{X} \ar@{^{(}->}[r] \ar[d]_{h} & X \ar[d]^{\overline{h}}\\
\mathcal{S} \ar@{^{(}->}[r]_{j_{\overline{\mathcal{S}}}} & \overline{\mathcal{S}}  
}
$$
\noindent o\`u $\overline{h}$ est projectif [Et 6, th\'eo (3.2.1)] ou [Et 5, chap I, th\'eo (3;3)].
En notant $E^i_{V} = R^i  f_{V{\textrm{rig}}^{\ast}}(X_{V}/\overline{\mathcal{S}}, j'^{\ast}(E))$, le th\'eor\`eme (3.1) prouve que $E^i_{V} \in F^a\mbox{-}\textrm{Isoc}^{\dag}(V/K)$.
On peut appliquer le [Et 7, (3.4.4)] ou [Et 5, chap II, (3.4.4)] qui fournit un isomorphisme\\

\noindent (3.2.2.2) $\qquad \qquad \qquad \qquad \qquad \widehat{E^i_{V}}\  \tilde{\longrightarrow}\  \mathcal{E}^i_{V}$\\

\noindent compatible aux Frobenius. Par le th\'eor\`eme 2 de [Et 4], les isomorphismes (3.2.2.1) et (3.2.2.2) assurent l'existence de $E^i \in F^a\mbox{-}\textrm{Isoc}^{\dag}(S/K)$ tel que

$$
\mathcal{E}^i = \widehat{E^i} \qquad \textrm{et}   \qquad E^i_{V} = j^{\ast}(E^i).
$$

\noindent L'unicit\'e de $E^i$ \`a isomorphisme pr\`es provient de la pleine fid\'elit\'e du foncteur d'oubli $F^a\mbox{-}\textrm{Isoc}^{\dag}(S/K) \rightarrow F^a\mbox{-}\textrm{Isoc}(S/K)$ \'etabli par Kedlaya [Ked 2, theo 1.1]. D'o\`u le th\'eor\`eme. $\square$\\

\vskip 3mm
\section*{4. Cas fini \'etale}

\noindent \textbf{Th\'eor\`eme (4.1)}. \textit{Soient $S$ un $k$-sch\'ema lisse et s\'epar\'e et $f : X \rightarrow S$ un $k$-morphisme fini \'etale. Alors, pour tout entier $i \geqslant 0$, $f$ induit des foncteurs canoniques }

$$
R^i  f_{\textrm{rig}\ast} : \textrm{Isoc}^{\dag}(X/K) \longrightarrow \textrm{Isoc}^{\dag}(S/K)
$$

$$
R^i  f_{\textrm{rig}\ast} : F^a\mbox{-}\textrm{Isoc}^{\dag}(X/K) \longrightarrow F^a\mbox{-}\textrm{Isoc}^{\dag}(S/K)
$$

\vskip 3 mm
\noindent \textit{et $\qquad \qquad  \qquad \qquad R^i  f_{\textrm{rig} \ast}(E)  = 0$ pour tout $i\geqslant 1$.} 

\vskip 3mm
\noindent \textit{D\'emonstration}. Soient $S_{0} = Spec\ A_{0} \hookrightarrow S$ un ouvert affine et $A_{1}$, $A_{2}$ deux $\mathcal{V}$-alg\`ebres lisses relevant $A_{0}$. On pose $S_{1} = Spec\ A_{1}$, $S_{2} = Spec\ A_{2}$ ; par la m\'ethode du [Et 6, th\'eo (3.1.1)] ou [ Et 5, chap I, th\'eo (3.1)] on a des compactifications $\overline{S}_{1} := P_{1}$, $\overline{S}_{2} := P_{2}$ de $S_{1}$ et $S_{2}$ et on note $\overline{S}_{0}$ l'adh\'erence sch\'ematique de $S_{0}$ plong\'e diagonalement dans $\overline{S}_{1} \times_{\mathcal{V}} \overline{S}_{2}$. En d\'esignant par $f_{0}$ la restriction de $f$ \`a $X_{0} = f^{-1}(S_{0})$ et par $\mathcal{S}_{1}$,  $\mathcal{S}_{2}$, $\overline{\mathcal{S}_{1}}$, $\overline{\mathcal{S}_{2}}$ les compl\'et\'es formels de $S_{1}$, $S_{2}$, $\overline{S}_{1}$, $\overline{S}_{2}$  respectivement, le th\'eor\`eme (3.1.1) de [Et 6] fournit des carr\'es cart\'esiens, $i = 1,2$,

$$
\xymatrix{
\mathcal{X}_{i} \ar[r] \ar[d]_{h_{i}} & \overline{\mathcal{X}_{i}} \ar[d]^{\overline{h}_{i}}\\
\mathcal{S}_{i} \ar@{^{(}->}[r]  & \overline{\mathcal{S}_{i}}
}
$$

\noindent o\`u $\overline{h}_{i}$ est fini, $h_{i}$ est fini \'etale et r\'el\`eve $f_{0}$ ; d'o\`u deux cubes commutatifs ($i = 1,2$)\\

$$
\xymatrix{
&\mathcal{X}_{i} \ar@{^{(}->}[rr] \ar@{.>}[dd]^(.4){h_{i}} |\hole && \overline{\mathcal{X}}_{i} \ar[dd]^{\overline{h}_{i}} \\
 \mathcal{X}_{1} \times_{\mathcal{V}}  \mathcal{X}_{2} \ar@{^{(}->}[rr] \ar[dd]_{h_{1} \times h_{2}} \ar[ur]^{u_{\mathcal{X}_i }} && \overline{\mathcal{X}}_{1} \times_{\mathcal{V}}  \overline{\mathcal{X}}_{2} \ar[dd] \ar[ur]_{u_{\overline{\mathcal{X}}_i }} \\
&\mathcal{S}_{i} \ar@{.>}[rr]^{\overline{h}_{1} \times \overline{h}_{2}}  && \overline{\mathcal{S}_{i}}  \\
 \mathcal{S}_{1} \times_{\mathcal{V}}  \mathcal{S}_{2} \ar[rr] \ar@{.>}[ur]_{u_{\mathcal{S}_i}} && \overline{\mathcal{S}}_{1} \times_{\mathcal{V}}  \overline{\mathcal{S}}_{2} \ar[ur]_{u_{\overline{\mathcal{S}}_i}}\ . \\
}
$$

\noindent Par le th\'eor\`eme [Et 7, (3.4.1)] ou [Et 5, chap II, (3.4.1)] on sait que pour $E \in \textrm{Isoc}^{\dag}(X/K)$ et $E_{0}$ sa restriction \`a $X_{0}$, alors $R^i f_{0 \textrm{rig}^{\ast}} (X_{0} / \overline{\mathcal{S}_{1}}, E_{0})$ et 
$R^i f_{0 \textrm{rig}^{\ast}}(X_{0} / \overline{\mathcal{S}_{2}}, E_{0})$ sont \'el\'ements de $\textrm{Isoc}^{\dag}(S_{0}/K)$ : de plus ils sont nuls pour $i \geqslant 1$ car $\overline{h}_{1}$ et $\overline{h}_{2}$ sont finis. De plus le th\'eor\`eme [Et 7, (3.4.4)] ou [Et 5, chap II, (3.4.4)] fournit des isomorphismes de changement de base

$$
u^{\ast}_{{\overline{\mathcal{S}}}_{i}} : f_{0 \textrm{rig} \ast} (X_{0} / \overline{\mathcal{S}_{i}}, E_{0}) \tilde{\longrightarrow} f_{0 \textrm{rig} \ast}(X_{0} / \overline{\mathcal{S}_{1}} \times_{\mathcal{V}} \overline{\mathcal{S}_{2}}, u^{\ast}_{\overline{\mathcal{X}}_{i}} E_{0}) ;
$$

\noindent d'o\`u un isomorphisme

$$
f_{0 \textrm{rig} \ast}(X_{0} / \overline{\mathcal{S}_{1}}, E_{0}) \tilde{\longrightarrow} f_{0 \textrm{rig}\ast}(X_{0} / \overline{\mathcal{S}_{2}}, E_{0}) \ ,
$$

\noindent et cet isomorphisme v\'erifie la condition de cocycles pour trois r\'el\`evements $S_{1}, S_{2}, S_{3}$ de $S_{0}$.\\

Par suite $f$ induit un foncteur

$$
f_{\textrm{rig}^{\ast}} : \textrm{Isoc}^{\dag}(X/K) \longrightarrow \textrm{Isoc}^{\dag}(S/K)
$$

\noindent puisque les constructions se recollent sur les ouverts de $S$. On pouvait aussi conclure en appliquant [Et 7, (3.4.8)] ou [Et 5, chap II, (3.4.8)].\\

La construction du Frobenius \'etant locale, on peut, pour montrer que $f_{\textrm{rig}^{\ast}}$ induit un foncteur 

$$
f_{\textrm{rig}^{\ast}} :  F^a\mbox{-}\textrm{Isoc}^{\dag}(X/K) \longrightarrow F^a\mbox{-}\textrm{Isoc}^{\dag}(S/K),
$$

\noindent supposer que $S$ est affine et lisse. La construction du th\'eor\`eme (2.1) s'applique ; le morphisme (2.1.17) est alors un isomorphisme car $F_{X/S}$ est un isomorphisme puisque $f$ est \'etale : l\`a on n'a pas besoin d'utiliser les r\'esultats de Ogus via le cas convergent (o\`u l'on avait suppos\'e $e \leqslant p-1$). On en d\'eduit directement que $\phi_{E_{i}}$ est un isomorphisme. D'o\`u le th\'eor\`eme. $\square$

\vskip 3mm

\noindent \textbf{Remarque (4.1.1)}. Tsuzuki a abord\'e dans [Tsu 1, theo (2.6.3)] la construction de $f_{\textrm{rig}^{\ast}}(X_{0}/\overline{\mathcal{S}},-)$ dans le cas fini \'etale, mais il n'\'etudie pas l'ind\'ependance par rapport \`a $\overline{\mathcal{S}}$ et ne prouve pas l'existence d'un $\mathcal{V}$- morphisme fini relevant le $f_{0}$ ci-dessus.

\vskip 3mm
\noindent \textbf{Th\'eor\`eme (4.2)}. \textit{Soient $S$ un $k$-sch\'ema s\'epar\'e de type fini, $E \in \textrm{Isoc}^{\dag}(S/K)$ et $f : X \rightarrow S$ un $k$-morphisme fini \'etale galoisien de groupe $G$}.\\

\noindent (4.2.1) \textit{Si $S$ est lisse sur $k$, alors, pour tout entier $i \geqslant 0$, on a des isomorphismes canoniques}\\

$\qquad \qquad H^i_{\textrm{rig}}(S/K, E)\  \tilde{\longrightarrow}\  (H^i_{\textrm{rig}}(S/K, f_{\textrm{rig}\ast} f^{\ast} E))^G$\\

$\qquad \qquad \qquad \qquad \qquad \tilde{\longrightarrow}\  (H^i_{\textrm{rig}}(X/K, f^{\ast} E))^G.$\\

\noindent (4.2.2) \textit{Si $k$ est parfait, ou si $S$ est affine et lisse sur $k$, alors, pour tout entier $i \geqslant 0$, on a des isomorphismes canoniques}\\

$\qquad \qquad H^i_{\textrm{rig},c}(S/K, E)\  \tilde{\longrightarrow}\  (H^i_{\textrm{rig},c}(S/K, f_{\textrm{rig}\ast} f^{\ast} E))^G$\\

$\qquad \qquad \qquad \qquad \qquad \tilde{\longrightarrow}\  (H^i_{\textrm{rig},c}(X/K, f^{\ast} E))^G.$\\

\noindent (4.2.3) \textit{Si $E \in F^a\mbox{-}\textrm{Isoc}^{\dag}(S/K)$ alors les isomorphismes de (4.2.1) et (4.2.2) sont compatibles \`a l'action du Frobenius.}

\vskip 3mm
\noindent \textit{D\'emonstration}. Par additivit\'e de la cohomologie rigide, avec ou sans supports, on peut supposer, pour le (1) et le (2), que $S$ est connexe.\\

\noindent \textit{Pour le (4.2.1)}, la suite spectrale de $\check{\mbox{C}}$ech en cohomologie rigide nous ram\`ene \`a $S$ affine et lisse sur $k$, $S = Spec\ A_{0}$. On choisit une $\mathcal{V}$-alg\`ebre lisse $A$ relevant $A_{0}$ et on reprend les notations utilis\'ees dans la preuve de (3.2.2): il existe un carr\'e cart\'esien de $\mathcal{V}$-sch\'emas formels

$$
\xymatrix{
\mathcal{X} \ar@{^{(}->}[r] \ar[d]_{h} & \overline{\mathcal{X}} \ar[d]^{\overline{h}} \\
\mathcal{S} \ar@{^{(}->}[r] & \overline{\mathcal{S}}
}
$$

\noindent et un syst\`eme fondamental $(V_{\lambda})_{\lambda} = (Spm\ A_{\lambda})_{\lambda}$ de voisinages stricts de $]S[_{\overline{\mathcal{S}}}$ dans $\overline{\mathcal{S}_{K}}$ et $\lambda_{0} > 1$ tel que pour $1 < \lambda \leqslant \lambda_{0}$ on ait un diagramme \`a carr\'es cart\'esiens

$$ \begin{array}{c}
\xymatrix{
\mathcal{X}_{K} \ar@{^{(}->}[r] \ar[d] _{h_{K}} & W_{\lambda} \ar@{^{(}->}[r] \ar[d]_{h_{\lambda}} & \overline{\mathcal{X}}_{K} \ar[d]^{\overline{h}_{K}}\\
\mathcal{S}_{K} \ar@{^{(}->}[r] & V_{\lambda} \ar@{^{(}->}[r] & \overline{\mathcal{S}}_{K}
}
\end{array}
\leqno{(4.2.1.1)}
$$

\noindent avec $\overline{h}_{K}$ fini, $h_{K}$ et $h_{\lambda}$ finis \'etales galoisiens de groupe $G$ [Et 7, (2.3.1)(2)] ou [Et 5, chap II, (2.3.1)(2)].\\

\noindent (4.2.1.2) Soit $E_{\lambda}$ un $\mathcal{O}_{V_{\lambda}}$-module localement libre de type fini. Pour $1 < \mu \leqslant \lambda$ on note

$$
\alpha_{\lambda_{\mu}} : V_{\mu} \hookrightarrow V_{\lambda} \quad , \quad \alpha_{\lambda} : V_{\lambda} \hookrightarrow \overline{\mathcal{S}}_{K}
$$

$$
\alpha'_{\lambda_{\mu}} : W_{\mu} \hookrightarrow W_{\lambda} \quad , \quad \alpha'_{\lambda} : W_{\lambda} \hookrightarrow \overline{\mathcal{X}}_{K},
$$

\noindent les immersions ouvertes et on pose

$$
j^{\dag}_{\lambda}\  E_{\lambda} = \displaystyle \mathop{\lim}_{\rightarrow} \alpha_{\lambda\mu^{\ast}}\  \alpha^{\ast}_{\lambda \mu}(E_{\lambda}),
$$

$$
j'^{\dag}_{\lambda}\  h^{\ast}_{\lambda}\  E_{\lambda} = \displaystyle \mathop{\lim}_{\rightarrow} \alpha'_{\lambda\mu^{\ast}}\  \alpha'^{\ast}_{\lambda \mu}\ h^{\ast}_{\lambda}\ (E_{\lambda}),
$$

$$
j^{\dag}\  E_{\lambda} =  \alpha_{\lambda^{\ast}}\  j'^{\dag}_{\lambda}\ E_{\lambda},
$$

$$
j'^{\dag}_{\lambda}\  h^{\ast}_{\lambda} (E_{\lambda}) = \alpha'_{\lambda^{\ast}}\  j'^{\dag}_{\lambda}\  h^{\ast}_{\lambda}\ (E_{\lambda}).
$$

\vskip 3mm
\noindent \textbf{Lemme (4.2.1.3)}. \textit{Avec les notations pr\'ec\'edentes on a des isomorphismes canoniques}

\begin{itemize}
\item[(i)] $(h_{\lambda^{\ast}}\  h^{\ast}_{\lambda}(E_{\lambda}))^G\  \tilde{\longrightarrow}\  E_{\lambda}.$
\item[(ii)] $(h_{\lambda \ast}\ h^{\ast}_{\lambda}\ j^{\dag}_{\lambda}\ E_{\lambda})^G \  \tilde{\longrightarrow}\  j^{\dag}_{\lambda}\ E_{\lambda}.$
\item[(iii)] $(\overline{h}_{K \ast}\ \overline{h}^{\ast}_{K}\ j^{\dag}\ E_{\lambda})^G\  \tilde{\longrightarrow}\  j^{\dag}\ E_{\lambda}.$
\end{itemize}

\vskip 3mm
\noindent \textit{D\'emonstration du lemme (4.2.1.3)}. 

\begin{itemize}
\item[(i)] Comme $E_{\lambda}$ est localement libre de type fini on a un isomorphisme
$$
h_{\lambda \ast}\ h^{\ast}_{\lambda}(E_{\lambda})\ \tilde{\longrightarrow}\  h_{\lambda \ast}\ h^{\ast}_{\lambda}(\mathcal{O}_{V_{\lambda}}) \otimes_{\mathcal{O}_{V_{\lambda}}} E_{\lambda},
$$
et l'action de $G$ sur le membre de gauche se fait par l'interm\'ediaire de $h_{\lambda \ast}\ h^{\ast}_{\lambda}(\mathcal{O}_{V_{\lambda}})$ puisque $G$ agit trivialement sur $E_{\lambda}$ : on est ramen\'e au cas $E = \mathcal{O}_{V_{\lambda}}$ qui a \'et\'e prouv\'e dans la proposition [Et 7, (2.3.1)] ou [Et 5, chap II, (2.3.1)].
\item[(ii)] On a des isomorphismes\\
$$
\begin{array}{r@{\ \, }c@{\ \, }l}
h_{\lambda \ast}\ h^{\ast}_{\lambda}\ j^{\dag}_{\lambda}\ E_{\lambda} & \simeq & h_{\lambda \ast}\ j'^{\dag}_{\lambda}\ h^{\ast}_{\lambda}\ E_{\lambda}\  [\mbox{B 3}, (2.1.4.7)]\\
&\simeq & j^{\dag}_{\lambda}\ h_{\lambda \ast}\ h^{\ast}_{\lambda}\ E_{\lambda}\  [\mbox{Et 7}, (3.1.4.1)] \ \mbox{ou} \ [\mbox{Et 5, chap II}, (3.1.4.1)]\\
&\simeq & h_{\lambda \ast}\ h^{\ast}_{\lambda}\ E_{\lambda}\ \otimes_{\mathcal{O}_{V_{\lambda}}} j^{\dag}_{\lambda}\ \mathcal{O}_{V_{\lambda}}\ [\mbox{B 3}, (2.1.3) (ii)].
\end{array}
$$
L'action de $G$ sur $h_{\lambda \ast}\ h^{\ast}_{\lambda}\ j^{\dag}_{\lambda}\ E_{\lambda}$ se fait par l'interm\'ediaire de $h_{\lambda \ast}\ h^{\ast}_{\lambda}\  E_{\lambda}$ puisque $G$ agit trivialement sur $j^{\dag}_{\lambda}(\mathcal{O}_{V_{\lambda}})$ : le (ii) r\'esulte alors du (i).
\item[(iii)] La preuve est semblable \`a celle du (ii) en utilisant cette fois [B 3, (2.1.4.8)] et [Et 7, (3.1.4.2)] ou [Et 5, chap II, (3.1.4.2)]. D'o\`u le lemme. $\square$
\end{itemize}

\vskip 3mm
Soit $E \in \textrm{Isoc}^{\dag}(S/K)$: on choisit le $V_{\lambda}$ comme ci-dessus de sorte qu'il existe un $\mathcal{O}_{V_{\lambda}}$-module localement libre et coh\'erent $E_{\lambda}$ tel que $j^{\dag} E_{\lambda}$ soit une r\'ealisation de $E$.\\

La cohomologie rigide $H^{\ast}_{\textrm{rig}}(S/K;E)$ est, pour $1 < \mu \leqslant \lambda$, la cohomologie des complexes

$$
\begin{array}{c}
\xymatrix{
 R \Gamma(V_{\mu};\  j^{\dag}_{\mu}  E_{\mu} \otimes_{\mathcal{O}_{V_{\mu}}} \Omega^{\bullet}_{V_{\mu}/K})\   \tilde{\leftarrow}\  R \Gamma (V_{\lambda};\  j^{\dag}_{\lambda}  E_{\lambda} \otimes_{\mathcal{O}_{V_{\lambda}}} \Omega^{\bullet}_{V_{\lambda}/K})\\
\qquad  \qquad \qquad \tilde{\rightarrow}\  M \otimes_{A^{\dag}_{K}} \Omega^{\bullet}_{A^{\dag}_{K}}\ , \\
}
\end{array} \leqno{(4.2.1.4)}
$$

\noindent o\`u la premi\`ere fl\`eche (resp. la deuxi\`eme) est un isomorphisme (resp. un quasi-isomorphisme), o\`u $M := \Gamma(V_{\lambda}; j^{\dag}_{\lambda}\ E_{\lambda})$ est un $A^{\dag}_{K}$-module projectif de type fini \`a connexion int\'egrable [B 3, (2.5.2)], o\`u $\Omega^1_{V_{\lambda}/K}$ est localement libre de type fini sur le faisceau coh\'erent d'anneaux $\mathcal{O}_{V_{\lambda}}$ [Et 7, (2.3.1)(2)] ou [Et 5, chap II, (2.3.1) (2)] et o\`u $\Omega^1_{A^{\dag}_{K}}$ est un $A^{\dag}_{K}$-module projectif de type fini [Et 4, 1.3].\\

De m\^eme la cohomologie rigide

$$
H^{\ast}_{\textrm{rig}}(S/K ; f_{\textrm{rig}\ast}\  f^{\ast}\ E)\quad\quad  (\textrm{resp.} H^{\ast}_{\textrm{rig}}(X/K ;  f^{\ast}\ E))
$$

\noindent est la cohomologie des complexes\\

\noindent (4.2.1.5)  $ \qquad \qquad   R\Gamma(V_{\lambda} ; h_{\lambda^{\ast}}\ h^{\ast}_{\lambda} (j^{\dag}_{\lambda}  E_{\lambda}) \otimes_{\mathcal{O}_{V_{\lambda}}} \Omega^{\bullet}_{V_{\lambda}/K})$\\

\noindent [resp. des complexes\\

\noindent (4.2.1.6)  $ \qquad \qquad   R\Gamma(W_{\lambda} ;  h^{\ast}_{\lambda} (j^{\dag}_{\lambda}  E_{\lambda}) \otimes_{\mathcal{O}_{W_{\lambda}}} \Omega^{\bullet}_{W_{\lambda}/K})\ ].$\\

\noindent Or la formule de projection, jointe au fait que $h_{\lambda}$ est \'etale, fournit des isomorphismes

$$
\begin{array}{c}
\xymatrix{
h_{\lambda^{\ast}} h^{\ast}_{\lambda} (j^{\dag}_{\lambda} E_{\lambda}) \otimes_{\mathcal{O}_{V_{\lambda}}} \Omega^{\bullet}_{V_{\lambda}/K} \simeq h_{\lambda^{\ast}} (h^{\ast}_{\lambda}\  j^{\dag}_{\lambda} E_{\lambda} \otimes_{\mathcal{O}_{W_{\lambda}}} h^{\ast}_{\lambda} (\Omega^{\bullet}_{V_{\lambda}/K})) \\
\qquad \qquad \qquad \qquad  \qquad \simeq\  h_{\lambda^{\ast}}(h^{\ast}_{\lambda}\  j^{\dag}_{\lambda}  E_{\lambda} \otimes_{\mathcal{O}_{W_{\lambda}}} \Omega^{\bullet}_{W_{\lambda}/K})\ ;\\
}
\end{array} \leqno{(4.2.1.7)}
$$

\noindent donc les complexes (4.2.1.5) et (4.2.1.6) sont quasi-isomorphes puisque $h_{\lambda}$ est fini.\\

Compte tenu du lemme (4.2.1.3) les isomorphismes

$$
H^i_{\textrm{rig}}(S/K;E)\  \tilde{\rightarrow}\  H^i_{\textrm{rig}}(S/K; f_{\textrm{rig}^{\ast}} f^{\ast} E)^G
$$

$\qquad \qquad \qquad \qquad  \qquad \qquad \quad \tilde{\rightarrow}\  H^i_{\textrm{rig}}(X/K; f^{\ast}E)^G$\\

\noindent s'\'etablissent comme [Et 1, (3.1.1)].\\

\noindent \textit{Pour le (4.2.2)}, comme la cohomologie rigide ne d\'epend que du sch\'ema r\'eduit sous-jacent, on supposera $S$ r\'eduit : si $k$ est parfait il existe alors un ouvert dense $U \hookrightarrow S$ avec $U$ affine et lisse sur $k$ et $Z := S \setminus U$ de dimension strictement plus petite que celle de $S$ [Et 2, d\'em. du th\'eo 3]. De plus $H^j(G, H^i _{\textrm{rig},c}(S/K, f_{\textrm{rig}\ast}\ f^\ast E)) = 0$ pour $j \geqslant 1$ [S\ 2, VIII, \S\ 2, cor 1 de prop 4] ; par fonctorialit\'e en $E$ de la cohomologie rigide \`a supports on en d\'eduit un morphisme de suites exactes longues

$$
\xymatrix{
\ar[r] & H^i_{\textrm{rig},c}(U, E_{\vert U}) \ar[d] \ar[r] & H^i_{\textrm{rig},c}(S,E) \ar[d] \ar[r] & H^i_{\textrm{rig},c} (Z, E_{\vert Z}) \ar[d] \ar[r] & \\
 \ar[r] & (H^i_{\textrm{rig},c} (U, f_{U\textrm{rig}\ast}f_{U}^{\ast}E_{\vert U}))^G \ar[r] & (H^i _{\textrm{rig},c} (S, f_{\textrm{rig}\ast}f^{\ast}E))^G \ar[r]  & (H^i_{\textrm{rig},c} (Z, f_{Z\textrm{rig}\ast}f_{Z}^{\ast}E_{\vert Z}))^G   \ar[r]  & 
}
$$
\\

\noindent Par r\'ecurrence sur la dimension on est ramen\'e \`a montrer l'isomorphisme du (4.2.2) pour $S$ affine et lisse sur $k$.\\

Reprenons les notations utilis\'ees pour la d\'emonstration du (4.2.1) et consid\'erons le diagramme commutatif \`a carr\'es cart\'esiens

$$
\begin{array}{c}
\xymatrix{
\mathcal{X}_{K} \ar@{^{(}->}[r] \ar[d]_{h_{K}} & W_{\lambda}  \ar[d]_{h_{\lambda}} & W_{\lambda} \setminus \mathcal{X}_{K} \ar[d]^{h'_{\lambda}} \ar@{_{(}->}[l]_{i'_{\lambda}} \\
\mathcal{S}_{K} \ar@{^{(}->}[r] & V_{\lambda}  & V_{\lambda} \setminus \mathcal{S}_{K} . \ar@{_{(}->}[l]_{i_{\lambda}} 
}
\end{array}
\leqno{(4.2.2.1)}
$$

\noindent La cohomologie \`a supports $H^{\ast}_{\textrm{rig},c}(S/K;E)$

$$
[\textrm{resp.} H^{\ast}_{\textrm{rig},c}(S/K; f_{\textrm{rig} \ast}\ f^{\ast}\ E), \textrm{resp.} H^{\ast}_{\textrm{rig},c}(X/K; f^{\ast}\ E)]
$$\\
\noindent est la cohomologie du complexe\\

\noindent (4.2.2.2) $\qquad  R\Gamma(V_{\lambda} ; j^{\dag}_{\lambda}  E_{\lambda} \otimes_{\mathcal{O}_{V_{\lambda}}} \Omega^{\bullet}_{V_{\lambda}/K} \longrightarrow i_{\lambda \ast} i^{\ast}_{\lambda}(j^{\dag}_{\lambda} E_{\lambda} \otimes_{\mathcal{O}_{V_{\lambda}}} \Omega^{\bullet}_{V_{\lambda}/K}))$\\

\noindent [resp.\\

\noindent (4.2.2.3) $R\Gamma(V_{\lambda} ; h_{\lambda \ast} h^{\ast}_{\lambda}(j^{\dag}_{\lambda}  E_{\lambda}) \otimes_{\mathcal{O}_{V_{\lambda}}} \Omega^{\bullet}_{V_{\lambda}/K} \rightarrow i_{\lambda \ast} i^{\ast}_{\lambda}(h_{\lambda \ast} h^{\ast}_{\lambda }(j^{\dag}_{\lambda} E_{\lambda})\otimes_{\mathcal{O}_{V_{\lambda}}} \Omega^{\bullet}_{V_{\lambda}/K}));$\\

\noindent resp.\\

\noindent (4.2.2.4) $R\Gamma(W_{\lambda} ;  h^{\ast}_{\lambda}(j^{\dag}_{\lambda}  E_{\lambda} \otimes_{\mathcal{O}_{W_{\lambda}}} \Omega^{\bullet}_{W_{\lambda}/K} \rightarrow i'_{\lambda \ast} i'^{\ast}_{\lambda}(h^{\ast}_{\lambda }(j^{\dag}_{\lambda} E_{\lambda})\otimes_{\mathcal{O}_{W_{\lambda}}} \Omega^{\bullet}_{W_{\lambda}/K}))].$\\

Le th\'eor\`eme de changement de base pour un morphisme propre [Et 7, th\'eo (3.3.2)] ou [Et 5, chap II, th\'eo (3.3.2)] fournit des isomorphismes

$$
 i_{\lambda \ast} i^{\ast}_{\lambda}(h_{\lambda \ast} h^{\ast}_{\lambda }(j^{\dag}_{\lambda} E_{\lambda})\otimes_{\mathcal{O}_{V_{\lambda}}} \Omega^{\bullet}_{V_{\lambda}}) \simeq  i_{\lambda \ast} i^{\ast}_{\lambda} h_{\lambda \ast} (h^{\ast}_{\lambda }(j^{\dag}_{\lambda} E_{\lambda})\otimes_{\mathcal{O}_{W_{\lambda}}} \Omega^{\bullet}_{W_{\lambda}})
$$
$$
\simeq  i_{\lambda \ast} h'_{\lambda \ast} i'^{\ast}_{\lambda} (h^{\ast}_{\lambda} j^{\dag}_{\lambda} E_{\lambda} \otimes \Omega^{\bullet}_{W_{\lambda}}) \simeq h_{\lambda \ast} i'_{\lambda \ast} i'^{\ast}_{\lambda} (h^{\ast}_{\lambda} j^{\dag}_{\lambda} E_{\lambda} \otimes_{\mathcal{O}_{W_{\lambda}}} \Omega^{\bullet}_{W_{\lambda}})\ ;
$$

\noindent donc via (4.2.1.7) les complexes (4.2.2.3) et (4.2.2.4) sont quasi-isomorphes puisque $h_{\lambda}$ est fini. L'isomorphisme (4.2.2) du th\'eor\`eme (4.2) en r\'esulte, compte tenu de (4.2.1.3) (ii).\\

\noindent \textit{Pour le (4.2.3)}, on peut supposer $S$ connexe affine et lisse sur $\mathcal{V}$ comme ci-dessus, dont on reprend les notations ainsi que celles de [Et 7, (2.3.1)(2)] (cf [Et 5, chap II, (2.3.1) (2)]). On fixe un rel\`evement $F_{A^{\dag}} : A^{\dag} \rightarrow A^{\dag}$
$$
\xymatrix{
\textrm{[resp.} F_{B^{\dag}} : B^{\dag} \ar[rr]_{1_{B^{\dag}}\otimes F_{A^{\dag}}}&& B^{\dag} \otimes A^{\dag}\ar[rr]^{\sim}_{F_{B^{\dag}/A^{\dag}}}&& B^{\dag}]\\
 }
 $$
 
 \noindent du Frobenius de $A_{0}$ [resp. de $B_{0}]$ comme dans [Et 5, (1.2)] : par extension des scalaires on en d\'eduit $F_{\hat{A}_{K}} : \hat{A}_{K} \rightarrow \hat{A}_{K}$ et $F_{\hat{B}_{K}} : \hat{B}_{K} \rightarrow \hat{B}_{K}$. On a vu en (1.2.4) qu'on dispose de carr\'es cart\'esiens o\`u $F_{\lambda \mu}$ et $F'_{\lambda \mu}$ sont finis :

$$
\begin{array}{c}
\xymatrix{
\mathcal{S}_{K} = Spm\  (\hat{A}_{K}) \ar@{^{(}->}[r] \ar[d]_{F_{\mathcal{S}_{K}}=Sp\ F_{\hat{A}_{K}}} 
& V_{\mu} = Spm(A_{\mu})  \ar[d]^{F_{\lambda \mu}}\\
\mathcal{S}_{K} = Spm\  (\hat{A}_{K}) \ar@{^{(}->}[r] & V_{\lambda} = Spm(A_{\lambda})
}
\end{array}
\leqno{(4.2.3.1)}
$$

\noindent et

$$
\begin{array}{c}
\xymatrix{
\mathcal{X}_{K} = Spm\  (\hat{B}_{K}) \ar@{^{(}->}[r] \ar[d]_{F_{\mathcal{X}_{K}}=Sp\ F_{\hat{B}_{K}}} 
& W_{\mu} = Spm(B_{\mu})  \ar[d]^{F'_{\lambda \mu}}\\
\mathcal{X}_{K} = Spm\  (\hat{B}_{K}) \ar@{^{(}->}[r] & V_{\mu} = Spm(B_{\lambda})\ ;
}
\end{array}
\leqno{(4.2.3.2)}
$$

\noindent plus pr\'ecis\'ement, \'etant donn\'e $\lambda$ on trouve $\mu$ de la fa\c{c}on suivante : en fixant des g\'en\'erateurs $\{ x_{i} \}$ de $B_{\lambda} $ sur $A_{\lambda}$ comme dans la preuve de [Et 7, (2.3.1) (2)], les \'el\'ements $F_{\hat{B}_{K}}(x_{i})$ sont entiers sur $A_{\lambda} \subset A^{\dag}_{K}$, donc a fortiori sur $B^{\dag}_{K} = \displaystyle \mathop{\lim}_{\rightarrow \atop{n}}\  B_{\lambda'}$ : il existe donc $\mu$, $1 < \mu \leqslant \lambda$ tel que pour tout $i$ on ait $F_{\hat{B}_{K}}(x_{i}) \in B_{\mu}$. Comme dans la preuve de [Et 7, (2.3.1) (2)] on peut aussi supposer que pour tout $i$ et tout $g \in G$ on a $g_{\hat{B}_{K}}(x_{i}) \in B_{\mu}$. Ainsi $F'_{\lambda \mu} : W_{\mu} \rightarrow W_{\lambda}$ (resp. $g_{\lambda} : W_{\lambda} \rightarrow W_{\lambda}$ est induit par $F_{B^{\dag}} : B^{\dag} \rightarrow B^{\dag}$ (resp. $g_{B^{\dag}} : B^{\dag} \rightarrow B^{\dag})$ ; pour prouver le lemme suivant:
\vskip 3mm

\noindent \textbf{Lemme (4.2.3.3)}. 

$$
g_{\lambda} \circ F'_{\lambda \mu} = F'_{\lambda \mu} \circ g_{\mu}.
$$
\noindent il suffit de prouver le 

\vskip 3mm

\noindent \textbf{Lemme (4.2.3.4)}. 

$$
g_{B^{\dag}} \circ F_{B^{\dag}} = F_{B^{\dag}} \circ g_{B^{\dag}}.
$$

\noindent Or $g \in G$ induit un morphisme $g_{X} : X \rightarrow X$ tel que $g_{X} \circ F_{X} = F_{X}
 \circ g_{X}$, puisque $g(x^q) = g(x)^q$ pour toute section $x$ de $O_{X}$ ; d'o\`u un diagramme commutatif\\
 
 $$
 \begin{array}{c}
 \xymatrix{
X  \ar@{}[drr] |{\rondI}  && X^{(q/S)}  \ar[ll]_{\pi_{X}} \ar@{}[drr] |{\rondII} && X  \ar[ll]_{F_{X/S}}^{\sim}  \ar@/_2pc/[llll]_{F_{X}} && \\
X \ar[u]^{g_{X}} && X^{(q/S)}   \ar[ll]^{\pi_{X}} \ar[u]^{g^{(q)}_{X}} && X \ar[ll]^{F_{X/S}}_{\sim}  \ar[u]_{g_{X}} &&\  .
}
\end{array}
\leqno{(4.2.3.5)}
$$

 \noindent Le carr\'e commutatif    $\ \rondI$ se rel\`eve en le carr\'e commutatif\\
 
 $$
\begin{array}{c}
\xymatrix{
B^{\dag} \ar[rr]^(.4){1_{B^{\dag}}  \otimes F_{A^{\dag}}} \ar[d]_{g_{B^{\dag}}} && B^{\dag} \otimes A^{\dag} \ar[d]^{g_{B^{\dag}} \otimes 1_{A^{\dag}}}\\
B^{\dag} \ar[rr]_(.4){1_{B^{\dag}}  \otimes F_{A^{\dag}}} && B^{\dag} \otimes A^{\dag} \ .
}
\end{array}
\leqno{(4.2.3.6)}
$$ 

\noindent Par l'\'equivalence de cat\'egories $B^{\dag} \longmapsto B_{0}$ de la cat\'egorie des 
$A^{\dag}$-alg\`ebres finies \'etales dans la cat\'egorie des $A_{0}$-alg\`ebres finies \'etales [Et 3, th\'eo 7], on rel\`eve le carr\'e commutatif     $\ \rondII$ en le carr\'e commutatif \\

 $$
\begin{array}{c}
\xymatrix{
B^{\dag} \otimes A^{\dag} \ar[d]_{g_{B^{\dag}} \otimes 1_{A^{\dag}}}\ar[rr]_{\sim}^(.6){F_{B^{\dag}/A^{\dag}}} && B^{\dag} \ar[d]^{g_{B^{\dag}}}&\\
B^{\dag} \otimes A^{\dag}  \ar[rr]^{\sim}_(.6){F_{B^{\dag}/A^{\dag}}} && B^{\dag} & .
}
\end{array}
\leqno{(4.2.3.7)}
$$ 

\noindent Par composition on a prouv\'e (4.2.3.4), donc (4.2.3.3).\\

Compte tenu de la commutation (4.2.3.3) et de la d\'efinition de la cohomologie rigide (resp. de la cohomologie rigide \`a supports compacts) donn\'ee en (4.2.1.4), (4.2.1.5), (4.2.1.6) [resp. en (4.2.2.2), (4.2.2.3), (4.2.2.4)] les isomorphismes (4.2.1) et (4.2.2) du th\'eor\`eme (4.2) sont compatibles \`a l'action du Frobenius. $\square$\\

\noindent Dans la preuve du th\'eor\`eme (4.2) on a montr\'e au passage :

\vskip 3mm
\noindent \textbf{Lemme (4.3)}. \textit{Si $S$ est un $k$-sch\'ema s\'epar\'e de type fini, $f : X \longrightarrow S$ est fini \'etale (non n\'ecessairement galoisien) et $E \in \textrm{Isoc}^{\dag}(X/K)$ on a des isomorphismes canoniques }

\begin{enumerate}
\item[(1)] $H^i_{\textrm{rig}}(X/K;E) \tilde{\longrightarrow} H^i_{\textrm{rig}}(S/K;f_{\textrm{rig}\ast}\ E)$.
\item[(2)] $H^i_{\textrm{rig},c}(X/K;E) \tilde{\longrightarrow} H^i_{\textrm{rig},c}(S/K;f_{\textrm{rig}\ast}\ E)$.
\item[(3)] \textit{Si de plus $E \in F^a\mbox{-}\textrm{Isoc}^{\dag}(X/K)$ les isomorphismes du (1) et (2) commutent \`a l'action de Frobenius}.
\end{enumerate}

\vskip 3mm
\noindent \textbf{Remarques (4.4)}.
\begin{itemize}
\item[(i)] Les r\'esultats du lemme (4.3) sont donn\'es par Tsuzuki dans [Tsu 1, cor (2.6.5) et (2.6.6)], sans pr\'ecisions de d\'emonstration, notamment pour le (2) du lemme : nous y avons utilis\'e le th\'eor\`eme de changement de base pour un morphisme propre [Et 7, (3.3.2)] ou [Et 5, chap II, (3.3.2)], qui n'est pas mentionn\'e  par Tsuzuki.
\item[(ii)] Le (4.2.2) du th\'eor\`eme (4.2) est une \'etape essentielle pour \'etablir la finiture de la cohomologie rigide \`a supports compacts \`a coefficients dans un $F$-isocristal surconvergent unit\'e \`a partir de la finitude de la cohomologie cristalline via la suite exacte longue de localisation en cohomologie rigide, la preuve de cet isomorphisme crucial n'appara\^it pas dans la d\'emonstration du th\'eor\`eme 6.1.2 de [Tsu 1].
\end{itemize}

\vskip10mm
\section*{5. Cas plongeable}

\textbf{5.1.} On suppose donn\'e un diagramme commutatif\\

$$
\xymatrix{
X \ar@{^{(}->}[r]^{j_{\mathcal{Y}}} \ar[d]_{f} & \mathcal{Y} \ar[d]^{\overline{h}} &\\
S \ar@{^{(}->}[r]_{j_{\mathcal{T}}} 		  & \mathcal{T} \ar[r]_{\rho} & Spf\  \mathcal{V}
}
$$

\noindent dans lequel $f$ est un morphisme de $k$-sch\'emas s\'epar\'es de type fini, $\overline{h}$ et $\rho$ sont des morphismes propres de $\mathcal{V}$-sch\'emas formels, $\overline{h}$ (resp. $\rho$) est lisse sur un voisinage de $X$ dans $\mathcal{Y}$ (resp. un voisinage de $S$ dans $\mathcal{T}$), $j_{\mathcal{Y}}$ et $j_{\mathcal{T}}$ sont des immersions. D\'esignons par $T$ (resp. $Y$) l'adh\'erence sch\'ematique de $S$ dans $\mathcal{T}$ (resp. de $X$ dans $\mathcal{Y}$), $\overline{f} : Y \rightarrow T$ le morphisme induit par $\overline{h}$, $i_{Y} : Y \hookrightarrow \mathcal{Y}$ l'immersion ferm\'ee, $X_{1} := \overline{f}^{-1}(S)$ et $f_{1} : X_{1} \rightarrow S$ le morphisme induit par $\overline{f}$.\\

On note $F_{S}$ (resp. $F_{X}$) le Frobenius de $S$ (resp. de $X$) (\'el\'evation \`a la puissance $q = p^a$ sur le faisceau structural) ; d'o\`u le diagramme commutatif

$$
\xymatrix{X \ar@/^1pc/[rrd]^{F_{X}} \ar@/_1,5pc/[rdd]_{f} \ar[rd]_{F_{X/S}} \\
&  X^{(q)} \ar[r]^{\pi_{X/S}} \ar[d]^{f^{(q)}} & X \ar[d]^{f}&\\
& S \ar[r]_{F_{S}} & S &.
}
$$

\vskip 3mm
\noindent \textbf{Th\'eor\`eme (5.2)}. \\
\textit{(5.2.1) Sous les hypoth\`eses (5.1) supposons que $\overline{h}^{-1}(S)= \overline{f}^{-1}(S) = X$; alors, pour tout entier $i \geqslant 0$, le morphisme $f$ induit un foncteur}

$$
R^i f_{\textrm{rig} \ast}(X/ \mathcal{T}, -) : F^a\mbox{-}\textrm{Isoc}^{\dag}(X/K) \longrightarrow F^a\mbox{-}\textrm{Isoc}^{\dag}(S/K).
$$

\noindent \textit{(5.2.2) Supposons donn\'es des morphismes}

$$
\xymatrix{
S' \ar@{^{(}->}[r]^{j'} & T' \ar@{^{(}->}[r]^{\rho'} & Spf \mathcal{V}
}
$$

\noindent \textit{o\`u $\rho'$ est un morphisme propre de $\mathcal{V}$-sch\'emas formels, $S'$ est un $k$-sch\'ema s\'epar\'e de type fini, $j'$ est une immersion et $\rho'$ est lisse sur un voisinage de $S'$ dans $\mathcal{T}'$. Alors le foncteur de (5.2.1) commute \`a tout changement de base s\'epar\'e de type fini $S' \rightarrow S$ : en particulier ce foncteur commute aux passages aux fibres en les points ferm\'es de $S$}.\\

\noindent \textit{D\'emonstration}.\\
\textit{Pour (5.2.1)}, soit $(E, \phi) \in F^a\mbox{-}\textrm{Isoc}^{\dag}(X/K)$; pour tout entier $i \geqslant 0$, on a d'apr\`es [Et 7, (3.4.4)] ou [Et 5, chap II, (3.4.4)] un isomorphisme

$$
F^{\ast}_{S} R^i f_{\textrm{rig} \ast} (X/ \mathcal{T}, E)  \displaystyle \mathop{\longrightarrow}^{\sim} R^if^{(q)}_{\textrm{rig} \ast}(X^{(q)}/\mathcal{T}, \pi^{\ast}_{X/S}(E)).
$$

\noindent L'identit\'e de $S$ induit un morphisme

$$
\xymatrix{
\theta^i : R^i f^{(q)}_{\textrm{rig} \ast}(X^{(q)}/ \mathcal{T}, \pi^{\ast}_{X/S}(E)) \ar[r] & R^i f_{\textrm{rig} \ast}(X/ \mathcal{T}, F^{\ast}_{X/S} \pi^{\ast}_{X/S}(E)) \ar@{=}[d]\\
&  R^i f_{\textrm{rig} \ast}(X/ \mathcal{T}, F^{\ast}_{X} (E)),
}
$$

\noindent et le Frobenius $\phi$ de $E$ induit un isomorphisme

$$
R^i f_{\textrm{rig} \ast}(X/ \mathcal{T}, F^{\ast}_{X} E) \displaystyle \mathop{\longrightarrow}^{\sim} R^i f_{\textrm{rig} \ast}(X/ \mathcal{T},  E).
$$

\noindent Par composition de ces trois morphismes on obtient le Frobenius de $R^i f_{\textrm{rig} \ast}(X/ \mathcal{T},  E)$\\

\noindent (5.2.3) $\qquad \qquad \qquad \phi^{i}: \  F^{\ast}_{S} R^i f_{\textrm{rig} \ast}(X/ \mathcal{T},  E) \longrightarrow R^i f_{\textrm{rig} \ast}(X/ \mathcal{T},  E) $\\

\noindent et il s'agit de prouver que $\phi^i$ est un isomorphisme : pour \c{c}a il suffit de prouver que c'est le cas pour $\theta^i$. On sait d\'ej\`a que $\theta^i$ est un morphisme d'isocristaux surconvergents  : d'apr\`es [B 3, (2.1.11) et (2.2.7)] il suffit de montrer que $\theta^i$ induit un isomorphisme dans la cat\'egorie convergente Isoc$(S/K)$ ; d'apr\`es [B-G-R, 9.4.3/3 et 9.4.2/7] il suffit de le v\'erifier apr\`es passage aux fibres de $\theta^i$ en les points ferm\'es $s$ de $S$. Pour un tel point $s$ notons $\mathcal{V}(s) = W(k(s)) \otimes_{W} \mathcal{V}$ et $K(s)$ le corps des fractions de $\mathcal{V}(s)$. D'apr\`es [Et 7, (3.4.4)] on a un diagramme commutatif

$$
\xymatrix{
R^i f^{(q)}_{\textrm{rig} \ast}(X^{(q)}/\mathcal{T}, \pi^{\ast}_{X/S}(E))_{s} \ar[r]^{\theta_{s}} \ar[d]^{\simeq} & R^i f_{\textrm{rig} \ast}(X/\mathcal{T}, F^{\ast}_{X}(E))_{s} \ar[d]_{\simeq}\\
R^i f_{s\ \textrm{rig} \ast}^{(q)} (X^{(q)}_{s}/\mathcal{V}(s), E_{X_{s}^{(q)}}) \ar@{.>}[r] \ar@{=}[d] & R^i f_{s\ \textrm{rig} \ast} (X_{s}/\mathcal{V}(s),F^{\ast}_{X_{s}} (E_{X_{s}})) \ar@{=}[d]\\
H^i_{\textrm{rig}} (X^{(q)}_{s}/K(s), E_{X^{(q)}_{s}}) \ar@{.>}[r]  \ar@{=}[d] & H^i_{\textrm{rig}}(X_{s}/K(s), F^{\ast}_{X_{s}}(E_{X_{s}})) \ar@{=}[d] \\
H^i_{\textrm{rig}, c} (X^{(q)}_{s}/K(s), E_{X^{(q)}_{s}}) \ar@{.>}[r]   & H^i_{\textrm{rig}, c}(X_{s}/K(s), F^{\ast}_{X_{s}}(E_{X_{s}})) 
}
$$

\noindent o\`u les fl\`eches verticales sont des isomorphismes ; or la fl\`eche horizontale inf\'erieure est un isomorphisme par [E-LS 1, prop 2.1, o\`u il faut supposer $X$ lisse sur $\mathbb{F}_{q}$ dans le cas de la cohomologie sans support]. D'o\`u (5.2.1).\\
L'assertion (5.2.2) r\'esulte de [Et 7, (3.4.4)]. $\square$

 \end{document}